\documentclass[12pt]{amsart}
\usepackage{amssymb,epsfig,amsmath,latexsym,amsthm,soul}

\usepackage{Bib-formating}

\usepackage[
colorlinks=true,
citecolor={black!80!black},
linkcolor={black!70!black},
urlcolor={black!80!black},
linktoc=section
]{hyperref}

\usepackage{graphicx}
\usepackage{pinlabel}
\usepackage{subcaption}
\usepackage{color}
\usepackage{tikz-cd}
\usepackage{color, xcolor, colortbl}
\theoremstyle{plain}
\newtheorem{theorem}{Theorem}[section]
\newtheorem{lemma}[theorem]{Lemma}

\newtheorem{proposition}[theorem]{Proposition}
\newtheorem{corollary}[theorem]{Corollary}

\theoremstyle{definition}
\newtheorem{definition}[theorem]{Definition}
\newtheorem{definition/construction}[theorem]{Definition/Construction}
\newtheorem{construction}[theorem]{Construction}
\newtheorem{remark}[theorem]{Remark}

\newtheorem{acknowledgements}[theorem]{Acknowledgements}

\newtheorem{data memo}[theorem]{Data Memo}

\newtheorem{notation}[theorem]{Notation}

\DeclareMathOperator{\Image}{Im}
\DeclareMathOperator{\Diff}{Diff}

\DeclareMathOperator{\Emb}{Emb}

\DeclareMathOperator{\id}{id}

\DeclareMathOperator{\tr}{tr}

\DeclareMathOperator{\LB}{LB}

\DeclareMathOperator{\I}{I}
\DeclareMathOperator{\C}{C}
\DeclareMathOperator{\grad}{grad}

\newcommand{\BR}{\mathbb R}

\newcommand{\mB}{\mathcal{B}}
\newcommand{\mC}{\mathcal{C}}
\newcommand{\mD}{\mathcal{D}}
\newcommand{\mE}{\mathcal{E}}
\newcommand{\mF}{\mathcal{F}}
\newcommand{\mG}{\mathcal{G}}

\newcommand{\mP}{\mathcal{P}}

\newcommand{\mR}{\mathcal{R}}

\newcommand{\mW}{\mathcal{W}}

\newcommand{\stwostwo}{S^2\times S^2}

\newcommand{\fF}{\mathfrak{F}}

\newcommand{\mfw}{(\mF,\mW)}

\begin{document}

\title{Sliding 2-handles over 3-handles and the barbell diffeomorphism}

\author{David Gabai}
\address{Department of Mathematics\\Princeton
University\\Princeton, NJ 08544}
\email{gabai@math.princeton.edu}

\author{Daniel Hartman}
\address{Department of Mathematics \\ Duke University \\ Durham, NC 27708}
\email{daniel.hartman@duke.edu}

\thanks{\today.
\newline\noindent\emph{Primary class:} 57R35
\newline\noindent\emph{secondary class:} 57R52, 57R50, 57N37
\newline\noindent\emph{keywords:} Smale, 4-sphere, diffeomorphism, pseudo-isotopy}

\begin{abstract} We give three fundamental examples highlighting the relationship between loops of spheres in $\#^k S^2\times S^2$, pseudo-isotopies, and diffeomorphisms of the 4-sphere. First, we describe an explicit pseudo-isotopy and Hatcher-Wagoner family that realizes the barbell diffeomorphism $\beta$ of the 4-sphere by constructing a finger-Whitney system for the family. From the Hatcher-Wagoner family, we are able to conclude that the example in \S 7 of ``Pseudo-isotopies of simply connected 4-manifolds'' generates $\beta$. Next, it is known that the square of $\beta$ is isotopic to id. We reprove this by building an explicit 2-parameter family of Morse functions that trivializes the Hatcher-Wagoner family for $\beta^2$. We do this by utilizing exchange points and families of handle slides. In the final example, we show that another Hatcher-Wagoner family, affectionately called the double cross, also generates $\beta$. \end{abstract}

\maketitle

\section{Introduction}
This paper focuses on the barbell diffeomorphism of $S^4$, and its connection with Hatcher-Wagoner families of Morse functions (see Definition \ref{D: HW families}). By ``the barbell diffeomorphism'', we really mean the isotopy class generated by the barbell diffeomorphism $\beta$ with supporting barbell shown in Figure \ref{F: 12BBGen} (Definitions \ref{D: Barbell} and \ref{D: BB map}). It turns out that this isotopy class contains several recent constructions of diffeomorphisms of $S^4$ including Watanabe's theta graph diffeomorphism \cite{wat18}, Gay's twists on unknotted Montesinos twins \cite{Gay, Kos2, Gay2}, and, as we will show in Corollary \ref{C: GGHKP}, Example 7.2 of \cite{GGHKP}. 

It is known that the barbell diffeomorphism has order 2 in $\pi_0\Diff(S^4)$, and that the square can be generated by a nontrivial loop of embedded circles in $S^1\times S^3$ \cite{BG1, GH}. In terms of Hatcher-Wagoner families, this means that a 1-eyed family of Morse functions and gradient-like vector fields with index 1 and 2 critical points generates $\beta^2 \in [\id]$ but that the eye cannot be directly canceled (see Lemma 3.6 of \cite{Laudenbach}). Theoretically, there must exist a relative homotopy from this family to one with no critical points since there are no nontrivial inertial pseudo-isotopies of $S^4$ \cite{Ce3}. To this end, we build a trivializing homotopy from a Hatcher-Wagoner family and show that it generates $[\beta^2]$. The family we construct consists of loops of $2$-spheres in $S^2\times S^2$, however it follows from the work in \cite{Gay, FGHK} that this can be ``traded'' for the nontrivial loop of circles. We leave this as an exercise for the interested reader. 

We also exhibit a new loop of spheres that realizes $[\beta]$. This loop, called the \textit{double cross}, consists of a loop of two 2-spheres in $S^2\times S^2 \# S^2\times S^2$, and is constructed in \S 5. In \cite{GGH1}, the authors, together with David Gay, constructed an invariant $\mD$ of loops of embedded spheres. The invariant $\mD$ is the $\mod 2$ sum of two quantities $\I$ and $\mC$ which count intersections between disks in a finger-Whitney system (Definition \ref{D: FW system}) for the loop. If we compare the two loops, namely the loop in \S3 and the double cross loop, using the work of \cite{GGH1}, we see that the first loop has $\I=1, \mC = 0$, while the double cross has $\I=0, \mC = 1$. Both have the same $\mD = \I+\mC$ invariant, but the individual quantities are different. To show that the double cross loop generates the barbell diffeomorphism, we again build a homotopy between the two Hatcher-Wagoner families. The homotopy used to interpolate between the double cross loop and the loop in \S 3 first eliminates one of the eyes in the double cross Hatcher-Wagoner family. What is interesting is that both of the homotopies constructed, the one in \S 4, and the one used in \S 5, realize different types of ``relations'' that can occur between the families. In the case where we trivialize $[\beta^2]$, we utilize exchange points (Definition \ref{D: exchange}) and families of handle slides, while in the second instance, we kill off one of the eyes and then homotope the resulting loop of spheres, using the machinery built in \cite{GGH1}. All of the phenomena shown in this paper are generic for 2-parameter families of Morse functions and associated gradient-like vector fields.  

\vspace{.4cm}
\noindent\textbf{AI statement:} Generative AI was not used in any part of this work.

\begin{acknowledgements}
We thank the organizers of the Georgia Topology
conference for the opportunity to present this work. Daniel Hartman would like to thank the Max Planck Institute for Mathematics where some of this work was conducted. David Gabai was partially supported by the National Science Foundation grants DMS-2304841 DMS-2603581. Daniel Hartman was supported by the National Science Foundation grant DMS-2602196.
\end{acknowledgements}

\section{Preliminaries}
Central to this paper are the relationships between $\pi_1 \Emb(\sqcup _kS^2, \#^k S^2\times S^2)$, $\pi_0\Diff(S^4)$ and $\pi_1(\fF,\mE)$, where $\fF$ denotes the space of smooth function from $S^4\times [0,1]$ to $[0,1]$ and $\mE$ the subset of those functions that have no critical points. We begin by reminding the reader of the general relationships that exist between the fundamental group of spaces of embeddings, pseudo-isotopies and 1-parameter families of Morse functions. We refer those interested in the details to \cite{HW}.

Let  
\[
\mP(X) = \{\phi\in \Diff(X\times [0,1]: \phi|_{X\times \{0\}} = (\id_X, 0) \}
\] denote the space of pseudo-isotopies of $X$, $\fF$ the space of smooth functions on $X\times [0,1]$ that have no critical points near $\partial X\times [0,1]$ and $\mE\subset \fF$ the subspace of $\fF$ consisting of smooth functions with no critical points. Now recall that $\mP(X)$ is homotopy equivalent to $\mE$, the two spaces being related by a fibration with contractible fiber \cite{Ce3, HW}, and that $\fF$ is itself contractible. Consequentially, we have that $\pi_{k+1}(\fF,\mE)\simeq \pi_k(\mP(X))$ for all $k\geq 0$. For our purposes, we will focus on the case when $k = 0$. Additionally, we will want gradient-like vector fields for these 1-parameter families. This additional information does not restrict us since the pair $(\fF,\mE)$ is homotopy equivalent to the pair $(\hat{\fF}, \hat{\mE})$ where $\hat{\fF}$ is the space of triples $(\xi, f, \mu)$ where $\xi$ is a gradient-like vector field for $f$ and $\mu$ is a metric on $X\times [0,1]$ \cite{HW}. We shall always implicitly work in $(\hat{\fF}, \hat{\mE})$ unless stated otherwise. As such, we will drop the ``hat'' and let $\fF$ and $\mE$ denote the space of triples.

The elements of $\mE$ are related to a pseudo-isotopy as follows. Given a triple $(\xi, f, \mu) \in \mE$, one gets a pseudo-isotopy $\phi$ by integrating the vector field $\xi/||\xi||$. Furthermore, if one starts with a pseudo-isotopy $\phi$, a metric $\mu$ on $X\times [0,1]$, one gets a triple $(\grad_\mu (\pi\circ \phi), \pi\circ \phi , \mu) $ where $\pi:X\times [0,1]\rightarrow [0,1]$ is the projection and $\grad_\mu (\pi\circ \phi)$ is the gradient of the function $\pi\circ \phi$. We will always implicitly assume there is a metric and will not refer to it and drop it from the notation. 

A consequence of $[\phi] = [\id_{X\times [0,1]}]$ is that the diffeomorphism $\phi$ is isotopic to an isotopy; a level preserving pseudo-isotopy. Furthermore, a pseudo isotopy is isotopic to $\id_{X\times [0,1]}$ if and only if any 1-parameter family $(q_t, \xi_t)$ in $(\fF,\mE)$ that interpolates between $\pi$ and $\pi\circ \phi$ and is homotopic relative to its endpoints to a path in $\mE$, i.e., $(q_t, \xi_t)$ can be deformed to a path with no critical points at all. In general, the question of if ``pseudo-isotopy implies isotopy'' for diffeomorphisms is false. There are some cases when it is known to be true: When $X$ is simply connected and $\dim(X)\neq 3,4$ \cite{Ce3, Ha1, Sm3}. The existence of an exotic diffeomorphism of a simply connected 4-manifold disproves the case when $\dim(X) = 4$ \cite{Rub98}, and an exotic diffeomorphism of the 4-sphere would disprove it when $\dim(X) = 3$ \cite{GGH2}.  

We will say a $1$-parameter family $(q_t, \xi_t)$ generates a pseudo-isotopy $\phi$ so long as $\xi_0 = \partial t$ and the flow of $\xi_1/||\xi_1||$ generates a pseudo-isotopy that is isotopic to $\phi$. Moreover, we can choose any representative of $[(q_i, \xi_i)]\in \pi_1(\fF,\mE)$ to generate $\phi$.

\begin{definition}\label{D: HW families}
    Let $(q_t, \xi_t)$ be a 1-parameter family with $(q_i, \xi_i)\in \mE$ for $i=0,1$. We say that the family is a \textit{Hatcher-Wagoner family} (HW family) with $k$-eyes if $q_t$ has only non-degenerate critical points or birth/death degenerate critical points such that
    \begin{enumerate}
        \item there are $k$ birth points and $k$ death points, and every birth point occurs before any death point.
        \item the index of any nondegenerate critical point is either 2 or 3, and
        \item  there is no gradient flow line of $\xi_t$ between any two nondegenerate critical points of the same index or to a birth/death point.
    \end{enumerate}
    A HW family is in \textit{nested eye position} if its Cerf graphic appears as shown in Figure \ref{F: HW family}.
\end{definition}
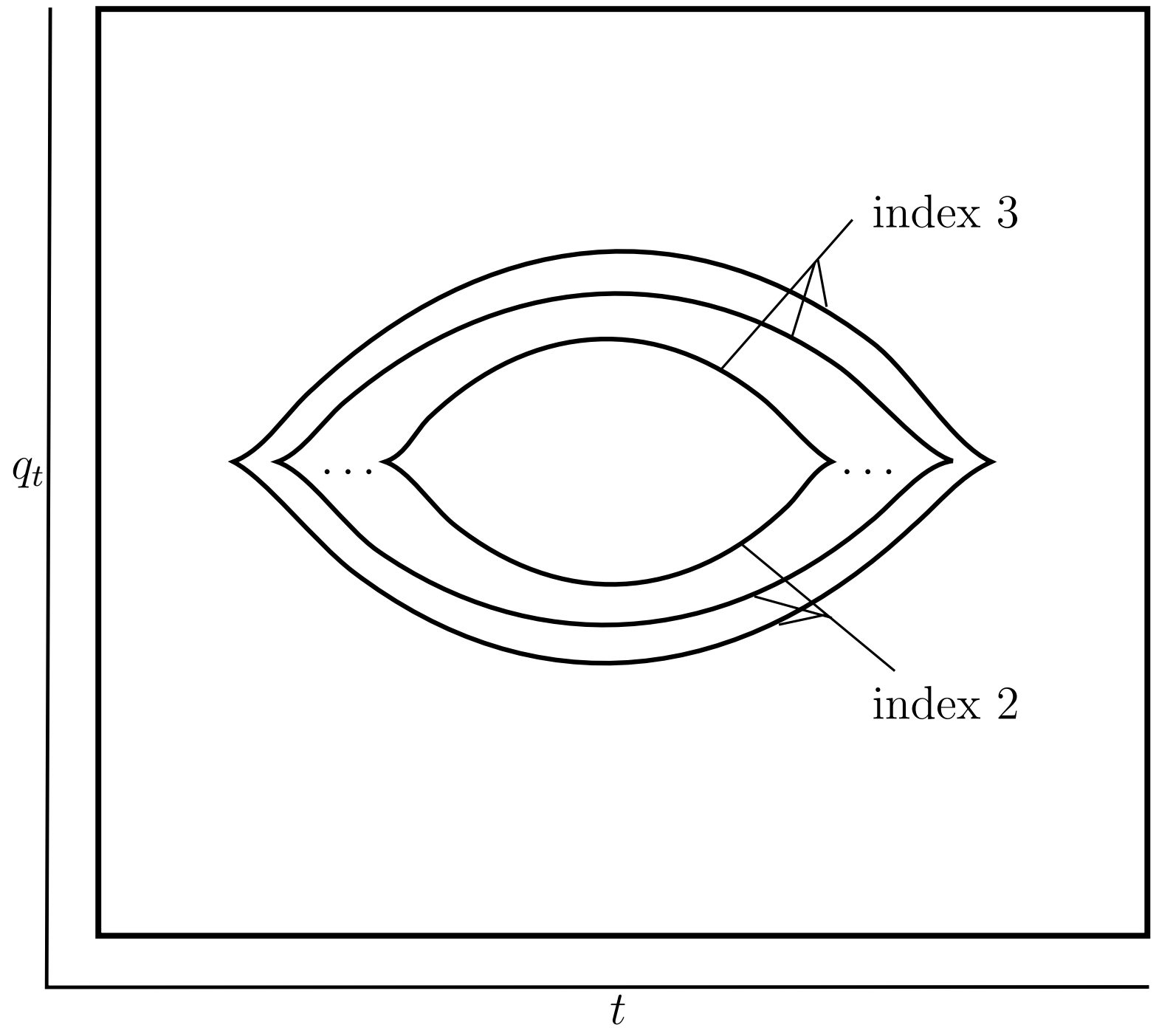
\begin{figure}
    \centering
    \labellist
    \pinlabel {$\cdots$} at 750 1260
    \pinlabel {$\cdots$} at 2020 1260
    \pinlabel {index $2$} at 2200 700
    \pinlabel {index $3$} at 2200 1900
    \pinlabel $q_t$ at -40 1260
    \pinlabel $t$ at 1400 -50
    \endlabellist
    \includegraphics[width=0.6\linewidth]{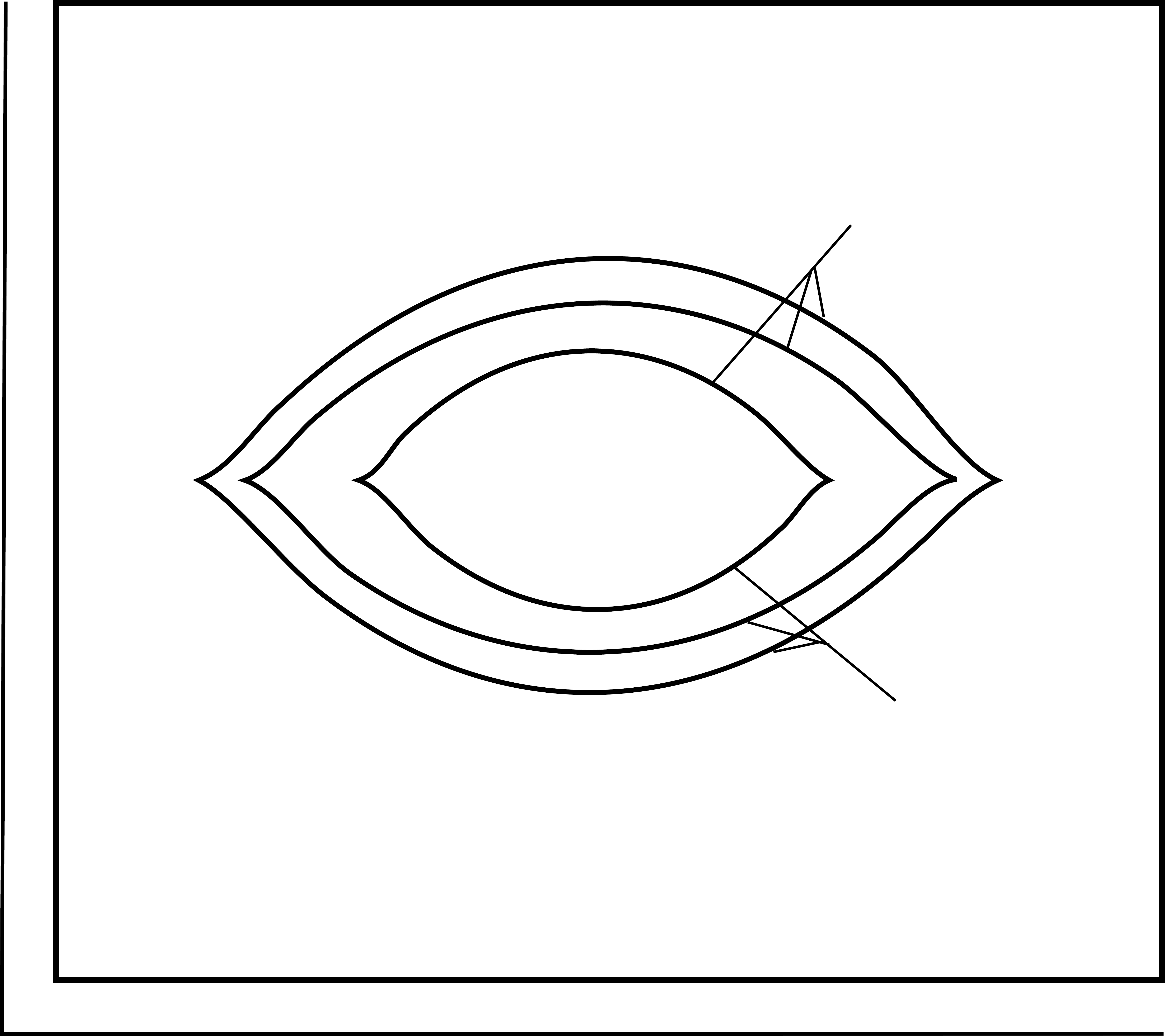}
    \vspace{.5cm}
    \caption{\label{F: HW family} Shown is a Cerf graphic for a nested eye Hatcher-Wagoner family.}
\end{figure}
It follows from \cite{HW} that when $X$ is simply connected and $\dim(X)\geq 4$, every pseudo-isotopy is generated by a $k$-eyed Hatcher-Wagoner family.
\begin{remark} 
A Hatcher-Wagoner family is always homotopic through HW familes to one in nested eye position \cite{HW}. 
\end{remark}

Given a nested eye HW family, we can modify $\xi_t$ so as to fix the descending spheres for the index 2 critical points in $X\times \{0\}$ to be constant except in an arbitarily small interval just before the death points, and that over this interval and until its death, there is a single flow line between any index 2 critical point and its canceling index 3 critical point. As such, we can modify the family by canceling the critical points over this interval, reducing to the case that the descending spheres for the index 2 critical points are constant. This gives us a canonical identificaiton of the level set above all the 2-handles with $X\#^k S^2\times S^2$. Furthermore, we can assume that this identifies the ascending spheres are given by $\mG = \sqcup_{i=1}^k (S^2\times \{r\})_i$.\footnote{Above the index 2-critical points, there is only one type of ascending sphere, namely the ones from the index 2 critical points. Similarly there is only one type of descending sphere.} The descending spheres for the index 3-critical points arrive into the level set as dual spheres for $\mG$. Near the birth/death points, they will be geometrically dual to $\mG$. Because of this, we can assume that near the birth and death points, the descending spheres are given by $\mR^{std} = \sqcup_{i=1}^k(\{g\}\times S^2)_i$. Now we can choose Morse coordinates for the index 3-critical points adapted to the descending manifolds in a consistent way so that flowing the characteristic embeddings down to the appropriate level set induces a map $\mR_t: \sqcup S^2 \rightarrow X\#^k S^2\times S^2$ for every $t$ in between the birth/death points \cite[Corollary 3.4]{Laudenbach}, \cite{milnorhcob}. As such we can associate to a nested eye HW family, a loop of embeddings $\mR_t$ based at $\mR^{std}$.

If given a loop $\mR_t$ based at $\mR^{std}$, we can construct a family of handlebodies and equip them with a family of Morse functions and gradient-like vector fields \cite{Gay, Ga3}. Then by identifying the handlebody with $X\times [0,1]$, we get a nested eye HW family on $X\times [0,1]$. 
\begin{proposition}\label{P: loop to HWfamily}
    Let $(q_t, \xi_t)$ be a nested eye HW family that induces a loop $\mR_t$ of spheres and let $(q'_t, \xi'_t)$ be the HW family given by building the handlebody using $\mG$ and $\mR_t$ and identifying this with $X\times [0,1]$. Then $(q_t, \xi_t)$ is relatively homotopic to $(q'_t, \xi'_t)$.
\end{proposition}
\begin{proof}
    We consider the concatenation of the two paths $(q_t, \xi_t)*(q'_{1-t}, \xi'_{1-t})$. After merging the death points of $(q_t, \xi_t)$ with the birth points of $(q'_t, \xi'_t)$, the middle level has a defining loop $\mR_t * \mR_{1-t}$ that is trivially isotopic to the constant loop based at $\mR^{std}$. As such, the family of gradient-like vector fields can be deformed so that the index 3 critical points and index 2 critical points are always in canceling position. Hence, all the eyes can be canceled, and therefore the composition is homotopic to a path in $\mE$. The result now follows from basic homotopy theory. 
\end{proof}

\begin{corollary}\label{C: homotopic loops}
    Homotopic loops of embeddings in $\Emb(\sqcup_k S^2, X\#^k S^2\times S^2)$ determine homotopic nested eye HW families.
\end{corollary}

\begin{remark}
    The above proposition says that the loop of spheres contains all the information about the nested eye HW family. 
\end{remark}

\begin{remark}
    Note that one actually extracts a loop of framed spheres from a HW family as the characteristic embeddings are framed embeddings. However since $\Omega^2 SO(2) \simeq *$, starting from a loop of spheres, one can always choose a framing for the loop, and this is unique up to isotopy of framed embeddings. We will often overlook the framings because of this. When needed for accuracy, we will denote a framed embedding by $\nu (\mR)$.
\end{remark}

Finally, let us review two different fibrations. First, recall that $$\Diff(X\#^k S^2\times S^2)\rightarrow \Emb_{\mR^{std}}^{fr}(\sqcup_k S^2, X\#^k S^2\times S^2)$$ is a locally trivial fiber bundle over its image with fiber $\Diff_\partial(X\#^k S^2\times S^2 \setminus \nu(\mR^{std}))$ \cite{Palais}. Now there is a natural inclusion of $\Diff_\partial(X\#^k S^2\times S^2 \setminus \nu(\mR^{std}))$ into $\Diff (X)$, 
\[
\centering
\begin{tikzcd}
    \Diff_\partial(X\#^k S^2\times S^2 \setminus \nu(\mR^{std}))\arrow[r, "\cup_{\id}"] &  \Diff(X).\\
\end{tikzcd} 
\]
\noindent The map is constructed by gluing $k$ copies of $S^1\times B^3$ along $\partial X\#^k S^2\times S^2 \setminus \mR^{std}$ and extending the diffeomorphism by $\id$ on each $S^1\times B^3$. This gives a well defined map from $\pi_0\Diff_\partial(X\#^k S^2\times S^2 \setminus \nu(\mR^{std}))$ to $\pi_0\Diff (X)$. If we precompose by the connecting homomorphism for the long exact sequence of the fibration,
\[
\delta: \pi_1\Emb(\sqcup_k S^2, X\#^k S^2\times S^2); \mR^{std}) \rightarrow \pi_0\Diff_\partial(X\#^k S^2\times S^2 \setminus \nu(\mR^{std})),
\]
then we get a map $\cup_{\id}\circ\delta:\pi_1\Emb(\sqcup_k S^2, X\#^k S^2\times S^2); \mR^{std}) \rightarrow \pi_0\Diff(X)$. This map is the \textit{parameterized surgery map}, and has appeared elsewhere in the literature \cite{Wall64, BG1, Ga3, Gay, Kos2, FGHK}.

The second fibration we will consider is 
\[
\begin{tikzcd}
    \mP(X)\arrow[r, "\rho"] & \Diff(X),
\end{tikzcd}
\]
\noindent where for a given pseudo-isotopy $\phi$, $\rho(\phi) = \phi|_{X\times \{1\}}$. Since $\pi_1(\fF, \mW) \simeq _{h.e.} \pi_0(\mP(X))$ and by Proposition \ref{P: loop to HWfamily} and Corollary \ref{C: homotopic loops} we have a well defined association of a homotopy class $[\mR_t] \in \pi_1 (\Emb(\sqcup_k S^2, X\#^k S^2\times S^2), \mR^{std})$ with $[(q_t, \xi_t)]\in \pi_1(\fF, \mW)$. Hence,

\begin{proposition}\label{P: ps map commutes}
    Let $\phi\in \mathcal{P(X)}$ be a pseudo-isotopy and let $\mR_t$ be a family of spheres extracted from a $k$-eyed HW family that generates $\phi$. Then $\cup_{\id}\circ \delta([\mR_t]) = [\rho\circ \phi]$ in $\pi_0\Diff(X)$.
\end{proposition}
\begin{proof}
    Let us denote by $X_k$ the manifold $X\#^k S^2\times S^2$. We start by applying the ambient isotopy extension theorem to the loop of spheres $\mR_t$. Let $\varphi_t$ be this extension. Note that this means $\mR_t = \varphi_t(\mR^{std})$. By appealing to the light bulb theorem \cite{Ga1}, we may assume that $\varphi_1(\mG) = \mG$ as well. From this data, we can build a pseudo-isotopy. We start with $X_k \times [0,1]\times [0,1]$ with coordinates labeled by $(x, t,y)$. Along the bottom $X_k\times [0,1]\times \{0\}$, we embed the track of $\mR_t$ and $\mG\times [0,1]\times \{0\}$. Along every vertical strip $X_k\times \{t\}\times [0,1]$, we apply the diffeomorphism $(\varphi_{t - ty}, y)$. On $X\times\{0\}\times [0,1]$ we have the identity and the track of $\varphi^{-1}$ on $X\times \{1\}\times [0,1]$. Moreover, along the top $X_k\times [0,1]\times \{1\}$ we have $\mR^{std}\times [0,1]\times \{1\}$. Now we can attach a family of 3-handles along the top of the square along $\mR^{std}\times [0,1]\times \{1\}$ and a family of upside down 2-handles along the bottom on $\mG\times [0,1]\times \{0\}$. The diffeomorphisms $\varphi_{t -ty}$ then fill in the ascending and descending manifolds over the square. At the ends, the handles are in canceling position with respect to each other. If we let $W_t = X_k\times \{t\}\times [0,1]\cup \text{handles}$, then $W_t$ is diffeomorphic to $X\times [0,1]$. Choose an identification $H: W_0\rightarrow X\times [0,1]$. Then we get an identification of $W_1$ by first using  $\Phi(x,y) =(\varphi_y(x), y)$ to identify $W_1$ to $W_0$, then composing with $H$. Thus we get a pseudo-isotopy of $X$ by $H\circ \Phi\circ H^{-1}$. It follows immediately from the construction that $H\circ \Phi\circ H^-1|_{X\times 1} = \cup_{\id}\circ\delta(\varphi_1)$. By Proposition \ref{P: loop to HWfamily}, the pseudo-isotopy produced is isotopic to the one that $\mR_t$ was originally extracted from, and hence $\cup_{\id}\circ \delta (\mR_t) = [\rho(\phi)]$.   
\end{proof}

\begin{notation}
    Moving forward, we will supress the notaiton $\cup_{\id}$ and just refer to $\delta$ as the parameterized surgery map.
\end{notation}

\section{1-parameter families of spheres for barbells.}

We begin by recalling some basic definitions.
\begin{definition}\label{D: Barbell}
    The \textit{model barbell} $\mB$ is the union of the two unit spheres in $\BR^3$ centered at $(-2,0,0)$ and $(2,0,0)$ and the arc from $(-1,0,0)$ to $(1,0,0)$ that runs along the x-axis. The 2-spheres are called the \textit{cuffs} and will be denoted by $C_1$ and $C_2$. The arc is referred to as the \textit{bar} and is denoted by $b$. The \textit{thickened model barbell} is the product $N_\epsilon(\mB) \times [-\epsilon, \epsilon]$ where $N_\epsilon(\mB)$ is an $\epsilon$-neighborhood of the model barbell in $\BR^3$. The \textit{mid-ball} $\Delta$ is the the subset $\{(0,y,z): y^2+z^2\leq \epsilon \}\times [-\epsilon, \epsilon]$ within the model barbell. We will abuse the notation and let $\mB: N_\epsilon(\mB)\rightarrow X$ denote an embedding of the thickened model barbell into a manifold $X$.
\end{definition}

\begin{remark}
    Note that an embedding of the thickened barbell manifold corresponds to an embedding of the model barbell such that both cuffs have trivial normal bundles. 
\end{remark}

\begin{definition}\label{D: BB map}
The thickened model barbell can also be defined as the complement of a neighborhood of two neatly embedded arcs in the 4-ball. By spinning the interior of one arc around the normal 2-sphere of the other, one gets an induced diffeomorphism via isotopy extension on the thickened model and this diffeomorphism restricts to the identity on a neighborhood of the boundary. We will let $\hat{\beta}$ denote a fixed representative of the isotopy class of diffeomorphisms that is generated by the spinning construction\footnote{Isotopy extension produces a diffeomorphism that is well defined up to the choice of extension and different extensions produce isotopic diffeomorphisms.}. See \cite{BG1} or \cite{FGHK} for a detailed explanation.

A \textit{barbell implantation} for an embedded barbell $\mB$ is the diffeomorphism defined as
    \begin{equation} \beta^\mB(x) = 
        \begin{cases}
            \mB\circ \hat{\beta}\circ \mB^{-1}(x) & x\in \Image(\mB)\\
            x & x\notin \Image(\mB)
        \end{cases}
    \end{equation}
    We say that $\mB$ is the \textit{supporting barbell} for $\beta^{\mB}$.
\end{definition}

\begin{figure}
    \centering
    \labellist
    \pinlabel $C_1$ at 250 -30
    \pinlabel $C_2$ at 900 -30 
    \endlabellist
    \includegraphics[width=0.5\linewidth]{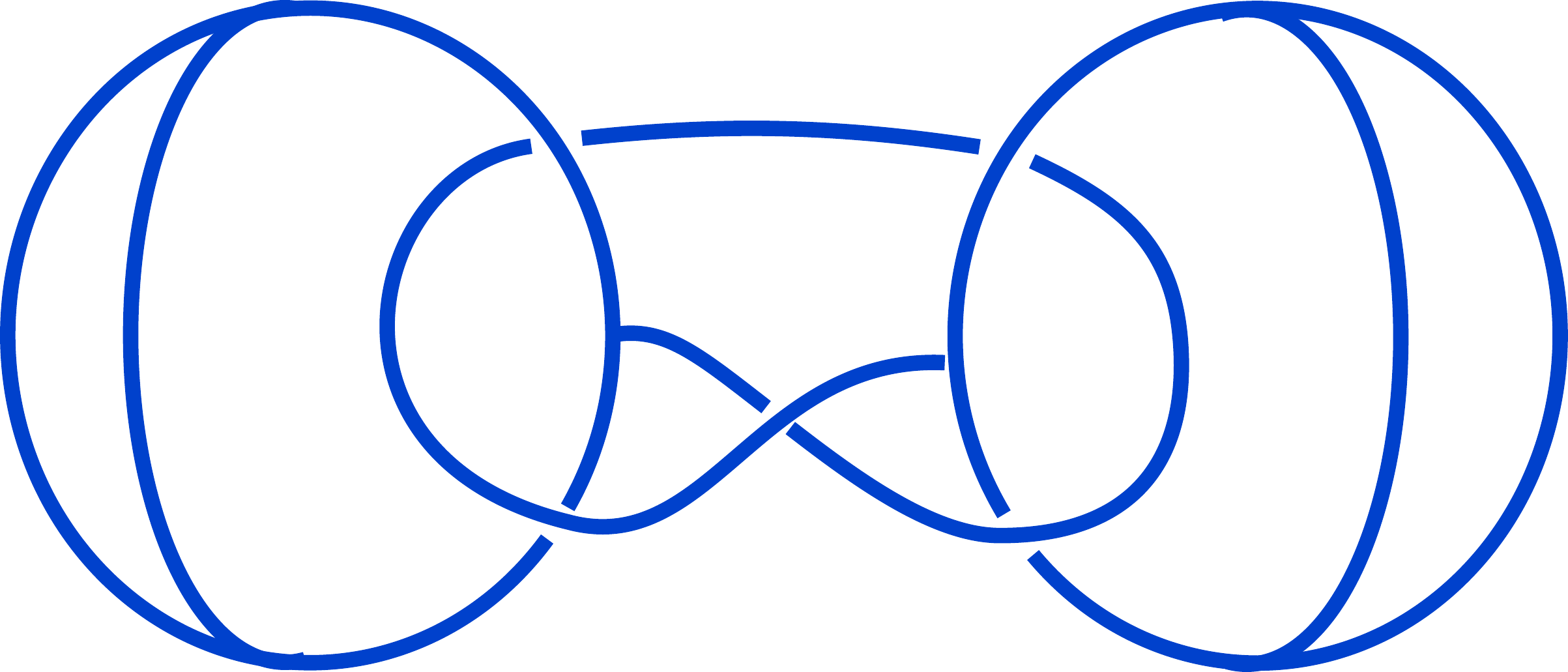}
    \vspace{.5cm}
    \caption{A Barbell embedding that gives an implantation that generates the 1-2 subgroup of $S^4$. \label{F: 12BBGen}}
\end{figure}

Our first objective is to construct a 1-parameter family of embedded spheres $\mR_t$ in $S^2\times S^2$ that generates the barbell diffeomorphism of $S^4$ for the barbell embeddeing shown in Figure \ref{F: 12BBGen}. In \cite{GH} it was shown that this diffeomorphism generates the 1-2 subgroup of $\pi_0\Diff^+(S^4)$ and that this subgroup contains all barbell diffeomorphisms with at least one cuff an unknotted $S^2$. 

\begin{remark}
    This subgroup is known as the \textit{1-2 subgroup} because every element arises from a Hatcher-Wagoner family with a
    single eye with index 1 and 2 critical points.
\end{remark}

\begin{notation}
    In the introduction, we let $\beta$ denote the barbell diffeomorphism supported on the embedding shown in Figure \ref{F: 12BBGen}. In this section, we will denote this diffeomorphism by $\beta^{1-2}$ when we refer to it and denote its embedding by $\mB^{1-2}$ to avoid confusion with the rest of the notation in this section. 
\end{notation}

\subsection{Constructing {$\mR_t$}}
We start by considering the barbell $\mB^{1-2}$ in $S^4\times \{1\}\subset S^4\times [0,1]$. To begin with, we choose a meridian $\mR^*$ for $C_1$ that bounds a disk $\mG^*$ that intersects $C_1$ in exactly 1 point. Next, we attach a 0-framed 5-dimensional 2-handle along $\mR^*$ and will work in the surgered manifold $S^4(\mR^*) \simeq S^2\times S^2$. We will let $\mR^{std}$ denote the ascending sphere for this 2-handle. Under the surgery, $\mG^*$ becomes a 2-sphere $\mG$ that is geometrically dual to $\mR^{std}$. The pair $\mR$ and $\mG$ allows us to identify $S^4(\mR^*)$ with $S^2\times S^2$ so that $\mR^{std} = \{g\}\times S^2$, $\mG = S^2\times \{r\}$. We fix this identification. 

\begin{remark}
    Note that we can recover $S^4$ from $S^2\times S^2$ by surgering $\mR^{std}$. If $\mR$ is any sphere in $S^2\times S^2$ isotopic to $\mR^{std}$, then surgering it produces a 4-manifold that is diffeomorphic to $S^4$ but not canonically.
\end{remark}

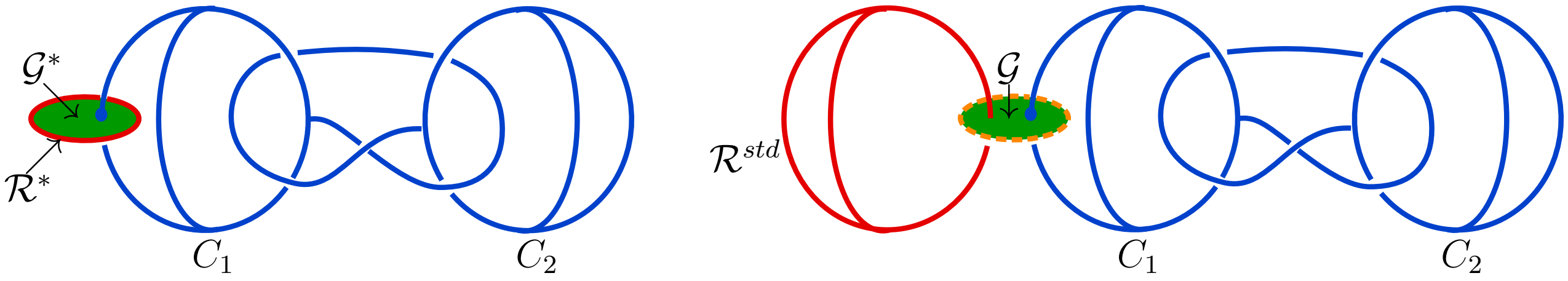
\begin{figure}
    \centering
    \labellist
     \hair 2pt
    \pinlabel $C_1$ at 400 -50
    \pinlabel $C_2$ at 1100 -50
    \pinlabel {$\nearrow$} [bl] at -20 120
    \pinlabel $\mR^*$ at 0 100
    \pinlabel $\mG^*$ at 30 360
    \pinlabel $\searrow$ at 70 290
    \pinlabel $\mR^{std}$ at 1550 170
    \pinlabel $\mG$ at 2120 360
    \pinlabel $\downarrow$ [t] at 2120 340
    \pinlabel $C_1$ at 2400 -50
    \pinlabel $C_2$ at 3100 -50
    \endlabellist
    \includegraphics[width=\linewidth]{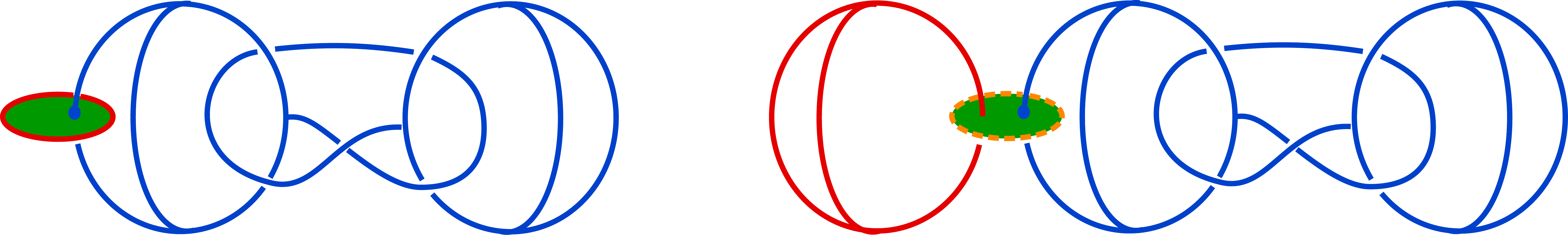}
    \vspace{.5cm}
    \caption{On the left we have depicted the barbell embedding $\mB^{1-2}$ in $S^4$ with $\mR^{*}$ and $\mG^*$ and cuffs $C_1$ and $C_2$ labeled. On the right, we have the barbell $\mB^{1-2}$ in $B^2\times S^2 \subset S^2\times S^2$, which is depicted as $S^4 \setminus \nu(S^1)$ (the dotted orange circle). We identify $B^2\times S^2$ with the subset of $S^2\times S^2$ that remains after removing $\{g'\}\times S^2$ so that $\mR$ is a \textit{parallel} sphere $\{g\}\times S^2$ for some $g \neq g'$. The green disk that is bounded by the dotted circle is what remains of $\mG$ after removing $\{g'\}\times S^2$, and $\mR^{std}$ is the red 2-sphere shown.
    \label{F: S4BBtoS2S2BB}}
\end{figure}

Since $\mR^*$ is contained in the complement of $\mB^{1-2}$, the barbell can be transported to $S^2\times S^2$. For convenience, we will draw $\mR^{std}$, $\mB^{1-2}$ in $B^2\times S^2$ that comes from removing a parallel $\{g'\}\times S^2$ that is disjoint from both. Removing this $\{g'\}\times S^2$ from $S^2\times S^2$ is equivalent to first surgering it to get $S^4$ and then removing the resulting $S^1\times B^3$ Thus, removing it allows us to conveniently view much of $S^2\times S^2$ in $S^4$ . See Figure \ref{F: S4BBtoS2S2BB}. 

\begin{figure}
    \centering
    \includegraphics[width=0.8\linewidth]{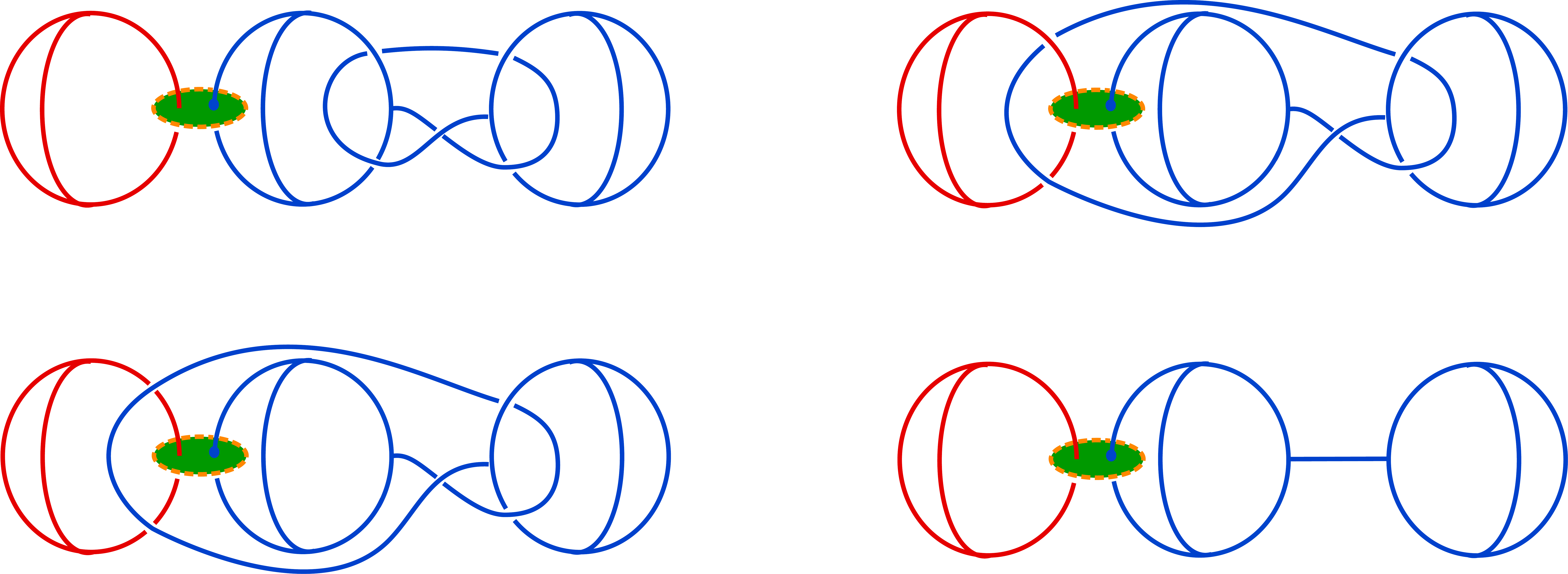}
    \caption{The following shows the isotopy of the barbell $\mB_t$ inside of $S^2\times S^2$.}
    \label{F: BBIsotopy}
\end{figure}

To generate the loop $\mR_t$, we will first define a path of barbell embeddings $\mB_t$ in $S^2\times S^2$ such that $\mB_0 = \mB^{1-2}$ and $\mB_1$ is in the complement of $\mR^{std}$ and induces the trivial diffeomorphism of $S^4$ after surgery. The path will be described in stages. First, we unlink the bar $b$ of $\mB^{1-2}$ with the cuff $C_1$ by using the hemisphere of $\mG$ not shown in the right figure of Figure \ref{F: S4BBtoS2S2BB} to isotope the bar until it no longer links $C_1$, but instead links $\mR^{std}$. This isotopy is performed in $S^2\times S^2$ and crosses $\{g'\}\times S^2$. The result is the top right of Figure \ref{F: BBIsotopy}. Next we unlink the bar from $\mR^{std}$ by pushing the bar through $\mR^{std}$. This is done is such a way that $\mR^{std}$ sweeps across the mid-ball for $\mB_t$. Then we isotope the resulting embedding until the bar does not link $C_2$ (see Figure \ref{F: BBIsotopy}).  

\begin{definition}
    We define $\beta_t := \beta^{\mB_t}$ as the induced path of diffeomorphisms given by applying the barbell map to the embedding $\mB_t$ at time $t\in [0,1]$.
\end{definition}
\begin{lemma}\label{L: BB_1 is trivial}
    $\beta_1$ is isotopic to the identity as a diffeomorphism of $S^2\times S^2$ and as an induced diffeomorphism of $S^4$. 
\end{lemma}
\begin{proof}
     Since $C_2$ bounds a 3-ball with an interior disjoint from $\mB_1$, the induced diffeomorphism $\beta_1$ is isotopic to the identity \cite[Lemma 9]{FGHK}. As shown in the bottom right of Figure \ref{F: BBIsotopy}, the obvious 3-ball is also disjoint from $\mR^{std}$. Therefore the isotopy trivializing $\beta_1$ in $S^2\times S^2$ fixes a neighborhood of $\mR^{std}$ and thus can be transported to $S^4$. 
\end{proof}

Given the path of diffeomorphisms $\beta_t$, we define the induced loop of spheres $\mR_t$ by the following formula:
\begin{equation}
    \mR_t = \beta_t(\mR^{std}).
\end{equation}

\begin{lemma}
    Let $\delta$ be the parameterized surgery map. Then $\delta([\mR_{1-t}]) = [\beta^{1-2}]$
\end{lemma}
\begin{proof}
    Since $\beta_t$ is an ambient isotopy extension of $\mR_t$, the lemma follows from  Propositions \ref{P: loop to HWfamily} and \ref{P: ps map commutes}.  
\end{proof}

\subsection{A finger-Whitney system for $\mR_t$.}
Now that we have the loop of spheres that generates $\beta^{1-2}$, we wish to extract the finger-Whitney system for this loop.

\begin{definition}\label{D: FW system}
    A finger-Whitney system is a tuple $(\mR, \mG, \mF, \mW)$ that consists of a pair of embeddings $\mR$ and $\mG$ in a 4-manifold $X^4$ such that both $\mR$ and $\mG$ have trivial normal bundles, together with two complete sets of framed Whitney disks $\mF$ and $\mW$ for $\mR\cup \mG$. 
\end{definition}

To construct the finger-Whitney system for $\mR_t$, we think of $\mG = S^2\times \{r\}$ as living in the codomain of $\beta_t$ and thus being stationary for the duration of $\beta_t$ and by definition, $\mR_t$. The support of $\mR_t$ is confined to the parameter values of $[0,1]$ for which the thickened bar $b$ pushes through $\mR^{std}$ in the domain of $\beta_t$ and by the construction of the path $\mB_t$, this is done so that $\mR^{std}$ only sweeps across the family of mid-balls $\Delta_t = \mB_t(\Delta)$. 

\begin{figure}
    \centering
    \includegraphics[width=0.8\linewidth]{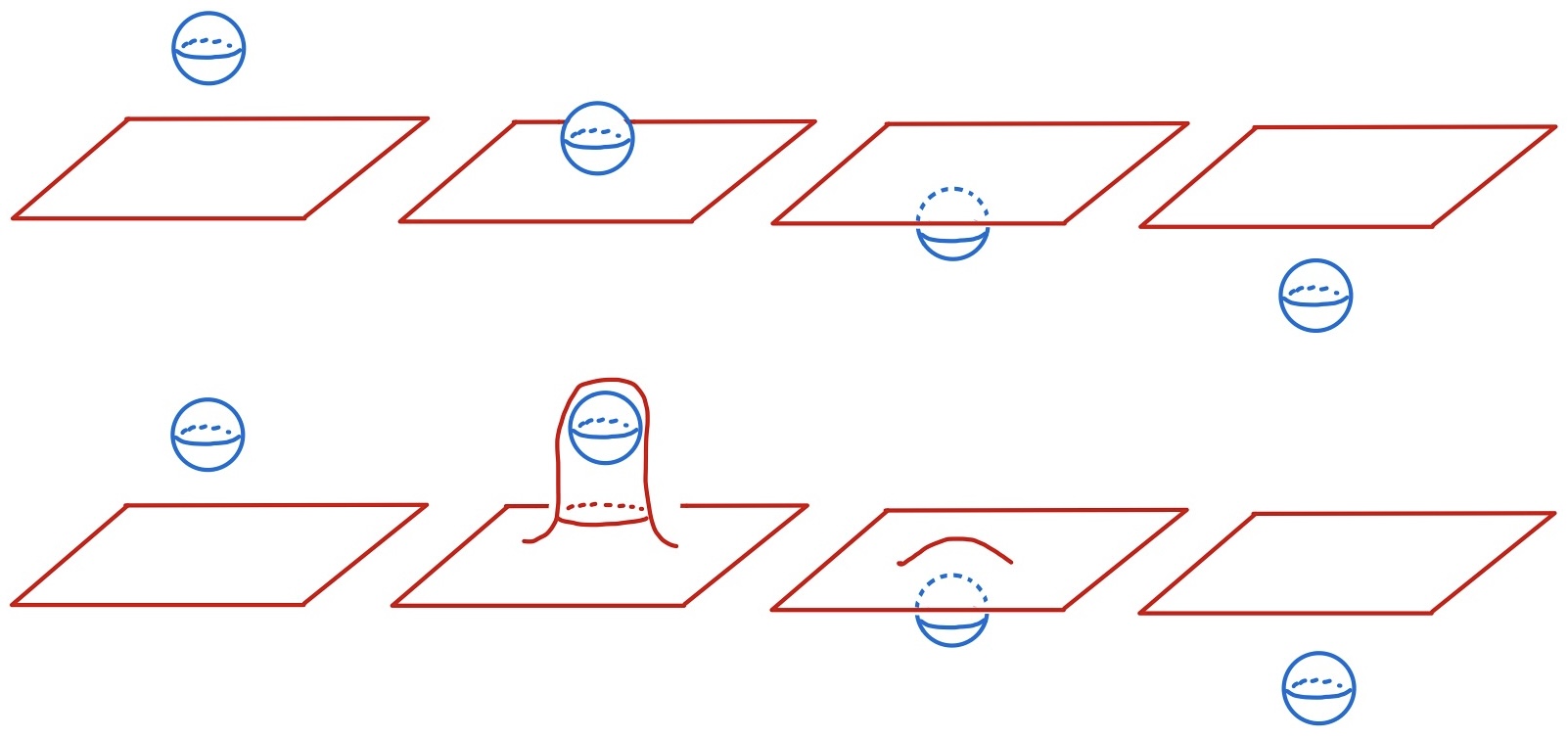}
    \caption{Shown is the local modification of the pair $\mB_t\cup \tilde{\mR}_t$ in the domain of $\beta_t$. Each frame is a local 3-dimensional slice in which $\tilde{\mR}_t$ appears as the horizontal red plane while $\mB_t$ intersects this 3-D slice in $\Delta_t$, shown as the blue sphere. The top row illustrates the original path while the bottom row illustrates the modified path.}
    \label{F: BarModify}
\end{figure}

Let $[t_0, t_1]$ be the subset of parameter values that goes from the top right to the bottom left in Figure \ref{F: BBIsotopy}. In the domain of $\beta_t$, the sphere $\mR^{std}$ can be viewed as a constant loop of spheres, while only $\mB_t$ moves. Let us denote by $\tilde{\mR}_t :=\mR^{std}$ the constant family in the domain of $\beta_t$, i.e., that $\beta_t(\tilde{\mR}_t)=\mR_t$. Now we modify the pair $\tilde{\mR}_t\cup \mB_t$ over the interval $[t_0, t_1]$ relative to $\{t_0, t_1\}$ so that from $[t_0, \frac{t_1+t_0}{2}]$ $\mB_t$ is constantly the embedding $\mB_{t_0}$ while $\tilde{\mR}_t$ protrudes a finger that swallows $\Delta_{t_0}$. Then from $[\frac{t_0+t_1}{2}, t_1]$ $\mB_t$ moves along its original path, while $\tilde{\mR}_t$ retracts back along the finger it protruded staying in the complement of $\mB_t$ for the duration. This modification can be done entirely locally and is shown in Figure \ref{F: BarModify}. We abuse notation and let $\tilde{\mR}_t$ denote the new modified loop and $\mR_t := \beta_t(\tilde{\mR}_t)$. 

\begin{remark}
Note that the deformation modifies the constant loop of $\mR^{std}$ to a non constant loop $\tilde{\mR}_t$. However, the deformation used is in fact a based homotopy of embeddings and thus $\tilde{\mR}_t$ is based homotopic to the constant loop. As such, the original loop $\mR_t = \beta_t(\mR^{std})$ is homotopic through embeddings to the new loop $\mR_t = \beta_t(\tilde{\mR}_t)$.  
\end{remark}

The main advantage of performing this deformation is freedom to choose how $\mR_t$ sweeps across the knotted mid-ball $\beta_t(\Delta_t)$ for $t\in [t_0, \frac{t_0+t_1}{2}]$. The ambiguity in this choice is controlled by $\Emb_\partial(B^2, B^3)$, which is contractible \cite{Ha1}, and therefore any two choices give homotopic loops of embeddings. Secondly, for $t\leq t_0$ and $t\geq\frac{t_0+t_1}{2}$, $\mR_t$ is geometrically dual to $\mG$. This means that the finger-Whitney moves all occur in $[t_0, \frac{t_0+t_1}{2}]$. 

\begin{figure}
    \centering
    \includegraphics[width=\linewidth]{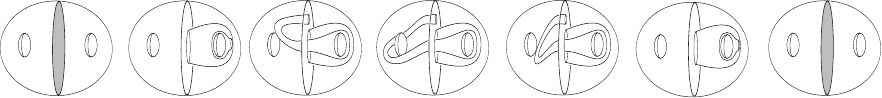}
    \caption{Shown is the 4-dimensional image of the mid-ball in the thickened model barbell $[-\epsilon, \epsilon]\times N_\epsilon (\mB)$ under the diffeomorphism $\hat{\beta}$. This is given by looking at the slices $\{x\}\times N_\epsilon(\mB)$ for several values of $x\in [-\epsilon, \epsilon]$.}
    \label{F: Midball}
\end{figure}

\begin{figure}
    \centering
    \labellist
    \small \hair 2pt
    \pinlabel $C_1$ at 300 800
    \pinlabel $T$ at 600 900
    \pinlabel $C_2$ at 2200 800
    \pinlabel $\{0\}\times B^2$ at 1450 150
    \endlabellist
    \includegraphics[width=0.5\linewidth]{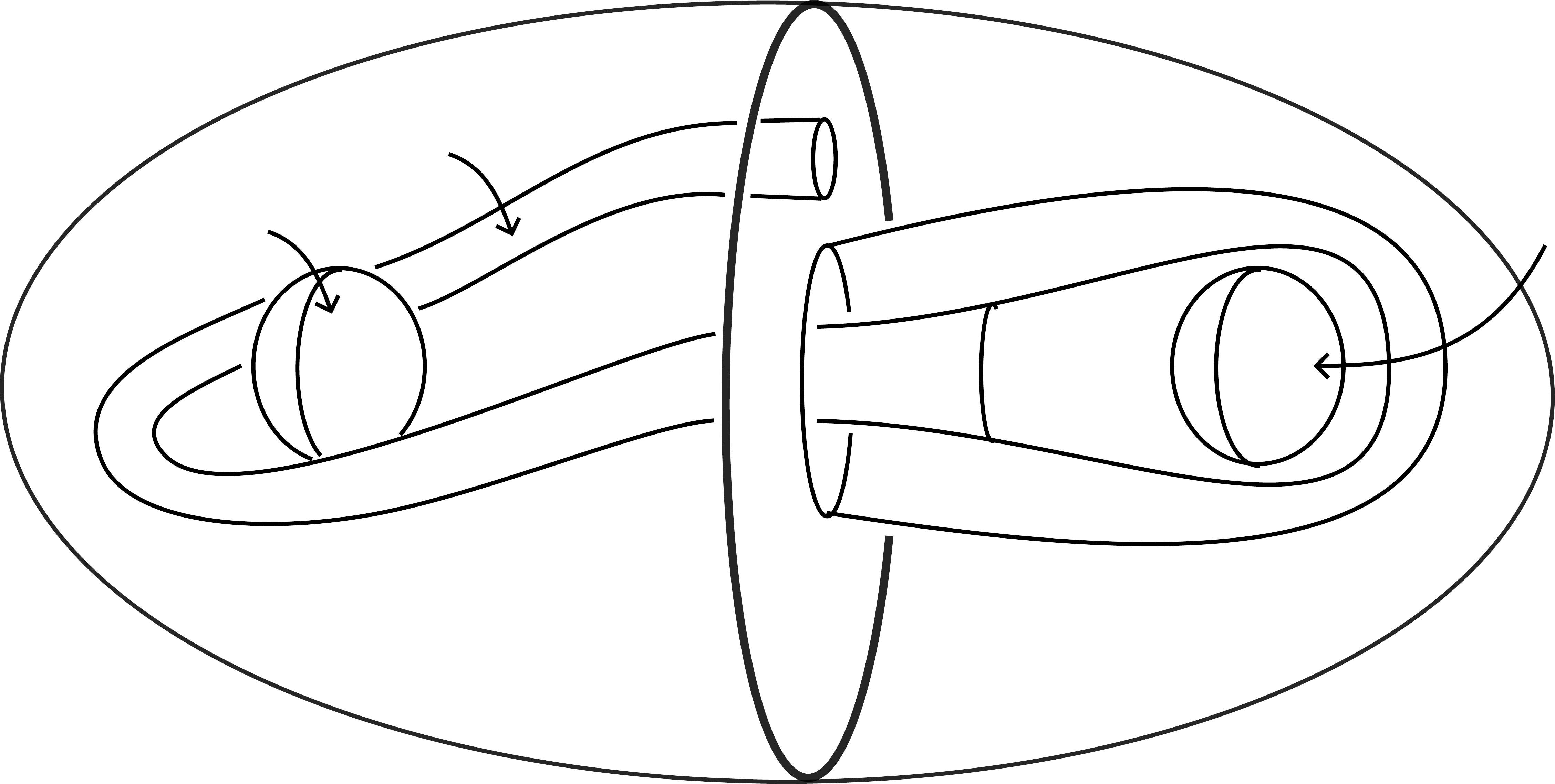}
    \caption{\label{F: midballhalfwaydisk} Shown is the slice of the mid-ball $\hat{\beta}(\Delta)$ that is contained in $\{0\}\times N_\epsilon (\mB)$. Note that in this slice, the disk $\hat{\beta(}\{0\}\times B^2)$ decomposes at the standard disk tubed to two copies of $C_2$. We have labeled the tube that swings around $C_1$ by $T$.}
    
\end{figure}
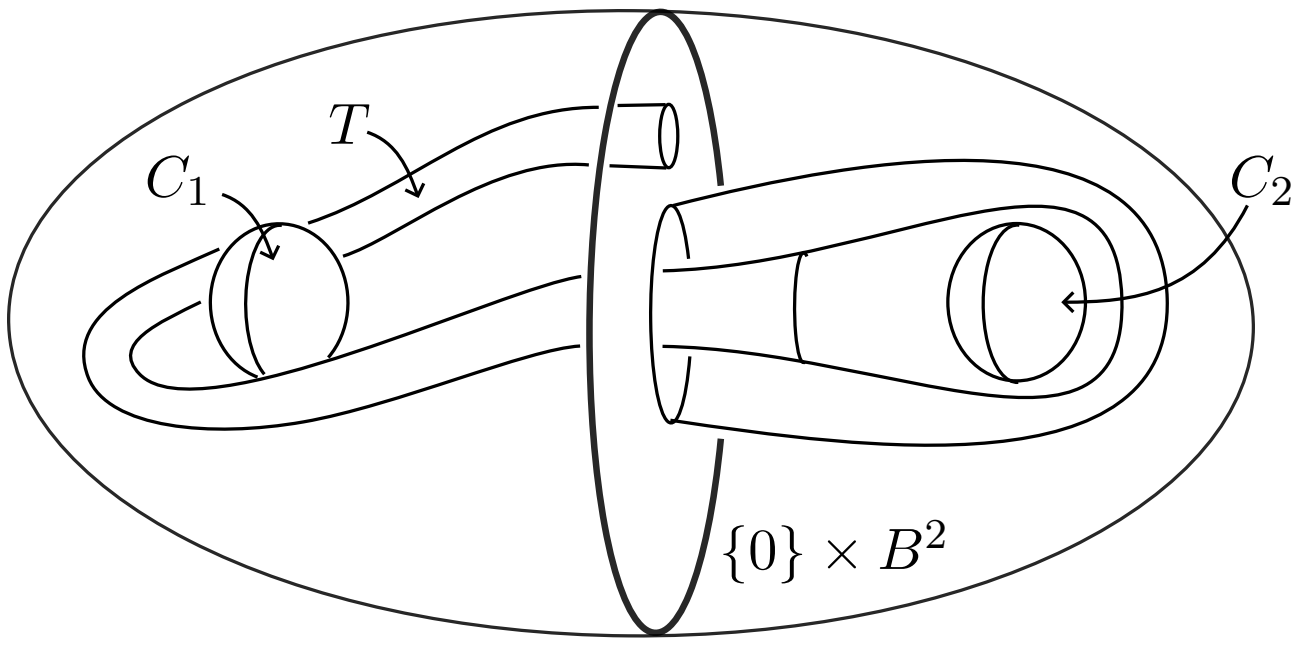

To get the finger-Whitney system, we recall from \cite{BG1} the structure of the mid-ball under the barbell diffeomorphism. This is shown in Figure \ref{F: Midball}. This figure depicts the thickened model barbell $[-\epsilon, \epsilon]\times N_\epsilon(\mB)$. Each frame is a slice $\{x\}\times N_\epsilon(\mB)$ for $x\in [-\epsilon, \epsilon]$. The mid-ball is shown as a family of disks $\{x\}\times B^2\subset \{x\}\times N_\epsilon(\mB)$. Importantly, the disk $\{0\}\times B^2$ is given by tubing the standard disk to two parallels of $C_2$. One of the tubes runs halfway along the bar $b$, while the other tube runs the full length of the bar from $C_2$ around $C_1$ and back half way to the disk. Furthermore, in $S^2\times S^2$ it can be arranged so that $\mG\cap\beta_{t_0}(\Delta) \subset \beta_{t_0}(\{0\}\times B^2)$ and that this intersection is confined to the tube that swings around $C_1$ \cite[Lemma 21]{FGHK}. We denote this tube by $T$ (see Figure \ref{F: midballhalfwaydisk}). 

\begin{figure}
    \centering
    \labellist
    \hair 2pt
    \pinlabel $\gamma$ at 550 4150
    \pinlabel $\searrow$ at 700 4100
    \pinlabel $\gamma$ at 3300 4450
    \pinlabel $\searrow$ at 3400 4350
    \pinlabel $f$ at 400 1850
    \pinlabel $\nwarrow$ at 300 1970
    \pinlabel $f$ at 4200 2400
    \pinlabel $\searrow$ at 4300 2300
    \pinlabel $w$ at 500 300
    \pinlabel $\nwarrow$ at 400 400
    \pinlabel $w$ at 4200 850
    \pinlabel $\searrow$ at 4300 750
    \endlabellist
    \includegraphics[width=\linewidth]{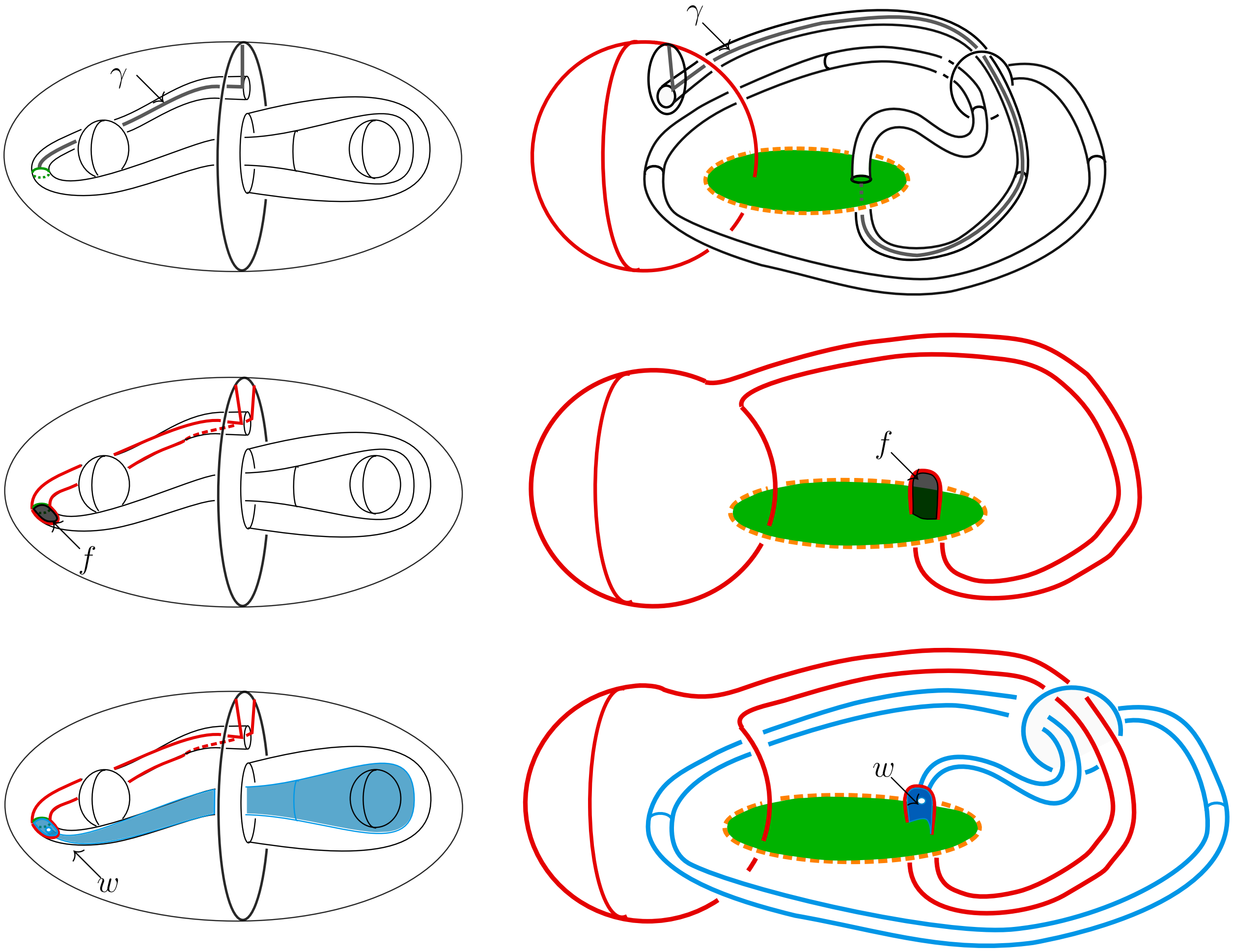}
    \caption{\label{F: FW data for 12 loop.} Shown is the data used to isotope $\mR$ across the subset $[-\delta, \delta]\times N_\epsilon (\mB)$ given in Construction \ref{Construction: FW Generator}. The images on the left shows the data in $\{0\}\times N_\epsilon (\mB)$ while on the right shows the image of the data under $\mB_{t_0}$ in $S^2\times S^2$. }
\end{figure}

\begin{construction}\label{Construction: FW Generator}
    We have arranged so that for $t\in[t_0, \frac{t_0+t_1}{2}]$ the path $\mR_t$ moves by is given by using a guiding arc to push a small finger of $\mR^{std}$ to $\beta_{t_0}(\Delta)$, then to sweep across the knotted 3-ball. We break the process of crossing the 3-ball into two stages. The first stage is given by pushing over the subset of a small thickening of $\beta_{t_0}(\{0\}\times B^2)$ containing the intersection with $\mG$, then pushing over the rest. The initial subset is constructed so that the isotopy given by sweeping across it is exactly the isotopy of a finger move followed by a Whitney move. 

\begin{figure}
    \centering
    \labellist
    \hair 2pt
    \pinlabel $U$ at 1300 1300
    \pinlabel $\swarrow$ at 1200 1170
    \pinlabel $T$ at 500 400
    \pinlabel $U$ at 2850 1300
    \pinlabel $T$ at 3350 1300
    \pinlabel $\searrow$ at 3450 1250
    \pinlabel $\searrow$ at 2950 1200
    \pinlabel $\downarrow$ at 420 850
    \pinlabel $C_1$ at 400 950
    \pinlabel $C_1$ at 4450 400
    \pinlabel $C_2$ at 4600 1300
    \pinlabel $C_2$ at 1900 400
    \pinlabel $\nwarrow$ at 1800 550
    \endlabellist
    \includegraphics[width=\linewidth]{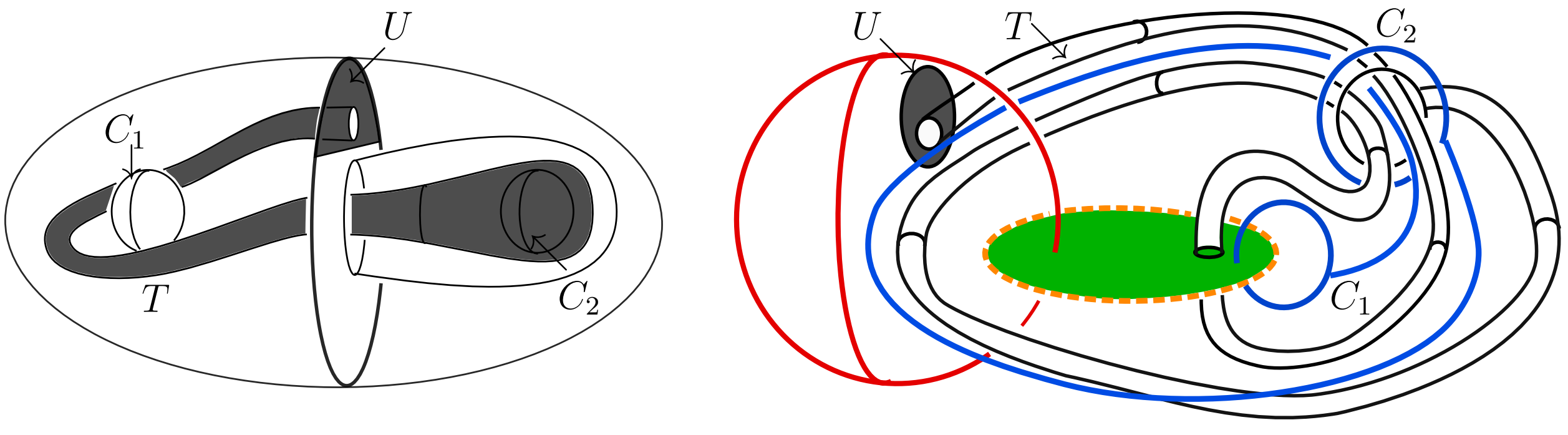}
    \caption{\label{F: subset isotopy} Shown is the subset $U\subset \hat{\beta}(\{0\}\times B^2)$ specified in Construction \ref{Construction: FW Generator} and its image under $\mB_{t_0}$ in $S^2\times S^2$.} 
\end{figure}
    
    We must first build this initial subset. Now, choose a small $\delta>0$ so that $\beta_{t_0}([-\delta, \delta]\times B^2)\subset \beta_{t_0}(\Delta)$ is a small thickening of the disk $\beta_{t_0}(\{0\}\times B^2)$. Recall from the paragraph above that $\hat{\beta}(\{0\}\times B^2)$ is constructed from a standard $\{0\}\times B^2$ by tubing it to $C_2$ twice, and that $T$ is the tube that sweeps around $C_1$ (see Figure \ref{F: midballhalfwaydisk}). We take the subset $U\subset \hat{\beta}(\{0\}\times B^2)$ as shown in the left picture in Figure \ref{F: subset isotopy}. The right picture in Figure \ref{F: subset isotopy} shows the image of $U$ in $S^2\times S^2$.
    
    The first isotopy is given by sweeping $\mR$ across $[-\delta, \delta]\times U$. We note that $U$ contains $T\cup C_2$ and that $[-\delta, \delta]\times U$ is also a 3-ball, so we are free to choose how we cross the ball as well. As such, we choose an arc $\gamma$ that run up the tube $T$ to $\mB_{t_0}^{-1}(\mG)$ and use this as the finger arc, pushing $\mR$ to $\mG\cap \beta_{t_0}(\{0\}\times B^2)$ and performing a finger move. The resulting finger disk $f$ can be seen as a subset of $T\cup C_1$. Then we use the rest of $T\cup C_1$ to perform a Whitney move to finish the isotopy. We have indicated this data in Figure \ref{F: FW data for 12 loop.}. The figures on the left side show the data in $\{0\}\times N_\epsilon(\mB)$ while the figures on the right shows the image of the data in $S^2\times S^2$.  \qed
\end{construction}

\begin{remark} 
Recall from \cite[\S 3.1]{GGHKP} that every framed Whitney disk has an associated sphere, called the \textit{Whitney sphere}, and that in \cite{GGH}, we reserved the name Whitney sphere for the spheres associated to the Whitney disks and called the spheres associated with the finger disks the \textit{finger spheres}. We keep with this convention here.
\end{remark}

\begin{proposition}\label{P: 12 FW system}
    The Whitney disk for the finger-Whitney system for $\mR_t$ is obtained from the finger disk by tubing to its own finger sphere.
\end{proposition}
\begin{proof}
    By Construction \ref{Construction: FW Generator}, the finger and Whitney disks are contained in $\beta_{t_0}(\{0\}\times U)$. The disks are formed after $\mR_t$ pushes along $\gamma$ to $\mG\cap \mB_{t_0}(\{0\}\times U)$. Let $L = \mG\cap \mB_{t_0}(\{0\}\times U)$ denote the intersection locus between $\mB_{t_0}$ and $\mG$. The interior of circle $L\subset \mB_{t_0}(\{0\}\times U)$ consists of the segment of $T$ that runs the length of the bar $b$ from $C_1$ to $C_2$ and the copy of $C_2$ with which it connects. Now, as $\mR_t$ crosses $[-\delta, \delta]\times U$, it will sweep across $L$. If we stop halfway through this isotopy, $\mR_t$ will intersect $\mG$ in two additional points, and the finger disk and Whitney disks consist of the closure of the complement of $\{0\}\times U \setminus \mR_t$ which is the interior of $L$. From this, it is clear that the Whitney disk is obtained from the finger disk by tubing $C_2$ to it using the rest of $T$. What remains to show is that $C_2$ is the finger sphere.
    
    To do this, we observe the embedding $\{0\}\times U$ in $S^2\times S^2$. This is illustrated in Figure \ref{F: FW data for 12 loop.}. Note that $\gamma$ runs along the bar $b$ from $C_2$ to $C_1$ and that the embedding $\mB_t$ embeds the bar so that this segment links $C_2$. It is clear that $C_2$ is a meridian 2-sphere for $\gamma$ and hence a finger-sphere.
\end{proof}

\begin{figure}
    \centering
    \includegraphics[width=\linewidth]{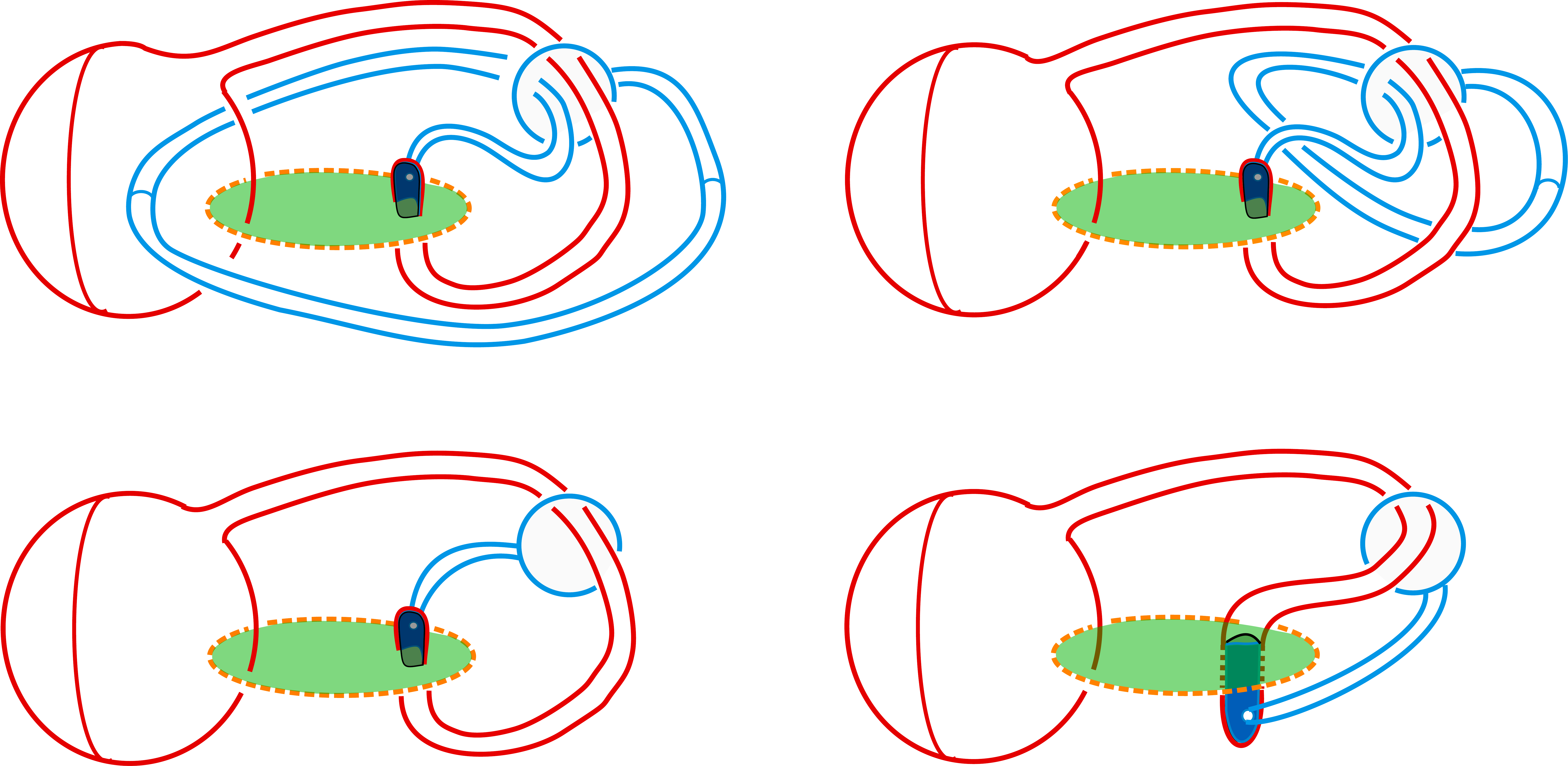}
    \caption{Shown is the isotopy of the finger-Whitney system for $\mR_t$ (top left) to one that agrees with Example 7.2 of \cite{GGHKP} (bottom right). The only difference here is the illustration of $S^2\times S^2$.}
    \label{F: FW isotopy to std}
\end{figure}

\begin{corollary}\label{C: GGHKP}
    The finger-Whitney system given in \S 7.2 of \cite{GGHKP} generates $[\beta^{1-2}]\in \pi_0\Diff(S^4)$. 
\end{corollary}
\begin{proof}
    First, let us note that Example 7.2 of \cite{GGHKP} has an equivalent representative finger-Whitney system where the Whitney disk and finger disk pair off the same two intersection points and the Whitney disk can be constructed by this finger disk by tubing to its finger sphere (see Example 0.9 of \cite{GGH}). Thus, the finger-Whitney system of $\mR_t$ and that of Example 7.2 differ in where the finger move is performed, and the tube used to construct the Whitney disk from the finger disk. Note that since $\pi_1(S^2\times S^2 \setminus \mR_0\cup\mG) = 1$ the finger arcs used to guide the finger moves can be made to agree. Since $\pi_1(S^2\times S^2 \setminus \mR_{t_*}\cup \mG)=1$ the tubes used to construct the Whitney disks can be made to agree after. Thus the two finger-Whitney systems are isotopic and therefore generate the same class $[\mR_t]\in \pi_1\Emb(S^2, S^2\times S^2)$ modulo $\pi_1 \LB$. Therefore they induce the same element $[\beta^{1-2}]\in \Diff(S^4)$.
    
    An isotopy of between the two systems is shown in Figure \ref{F: FW isotopy to std}. In going from the top left to the top right in the figure, the tube $T$ swings around $\mG$ in $S^2\times S^2$ and thus cuts across $\{g'\}\times S^2$ (not shown). Then $T$ continues to move, going from the top right to the bottom left. Finally, by rotating the finger move of $\mR$ into $\mG$ near $\mG$ gives the final isotopy going from the bottom left to the bottom right.
\end{proof}

\section{Barbells and exchange points.}
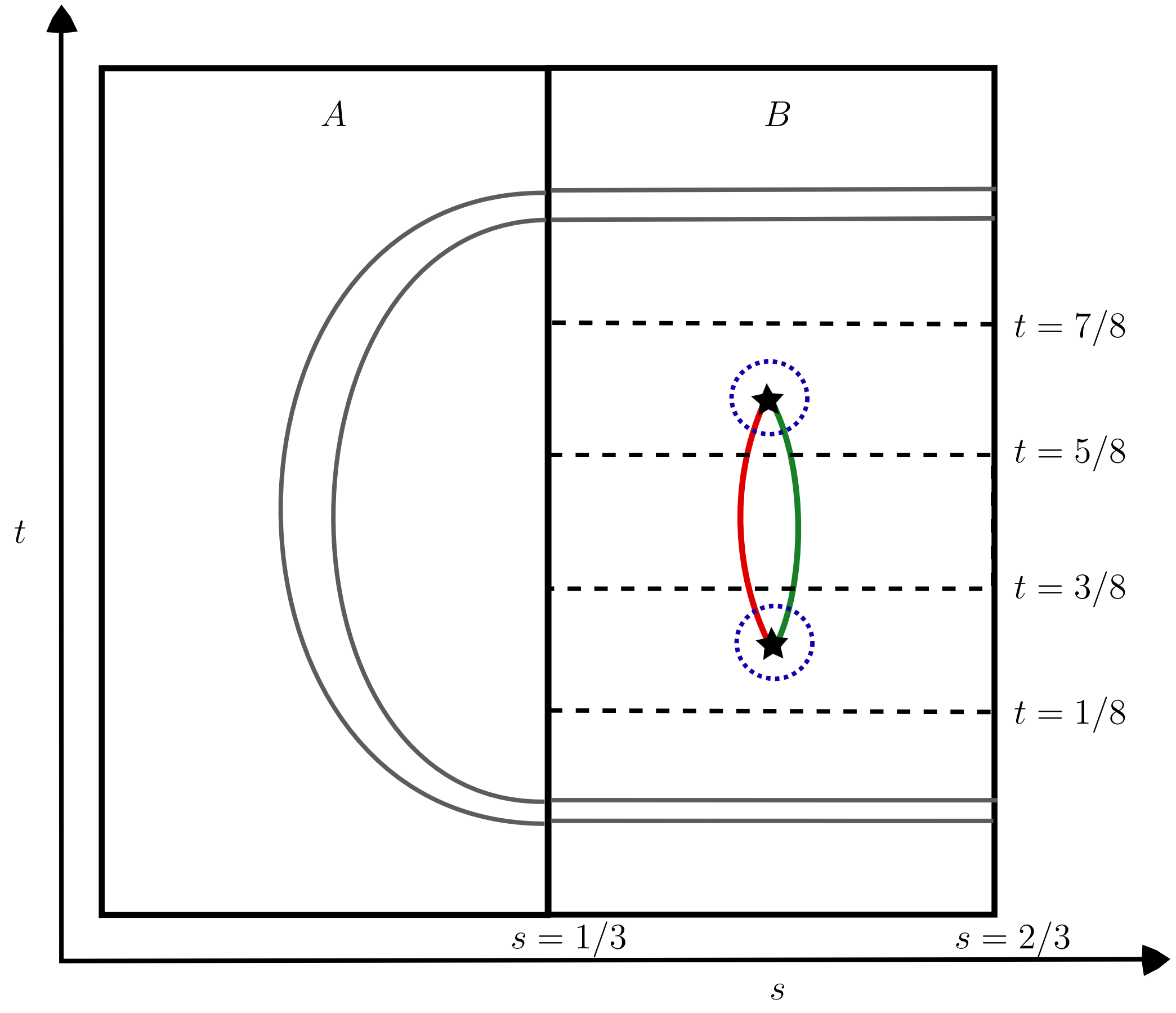
\begin{figure}
    \centering
    \labellist                             
    \pinlabel $A$ at 550 1650
    \pinlabel $B$ at 1400 1650
    \pinlabel $t$ at -50 850
    \pinlabel $s$ at 1400 -30
    \pinlabel {$t=1/8$} at 1960 500
    \pinlabel {$t=3/8$} at 1960 740
    \pinlabel {$t=5/8$} at 1960 1000
    \pinlabel {$t=7/8$} at 1960 1240
    \pinlabel {$s = 1/3$} at 1000 70
    \pinlabel {$s = 2/3$} at 1850 70
\endlabellist
    \includegraphics[width=0.8\linewidth]{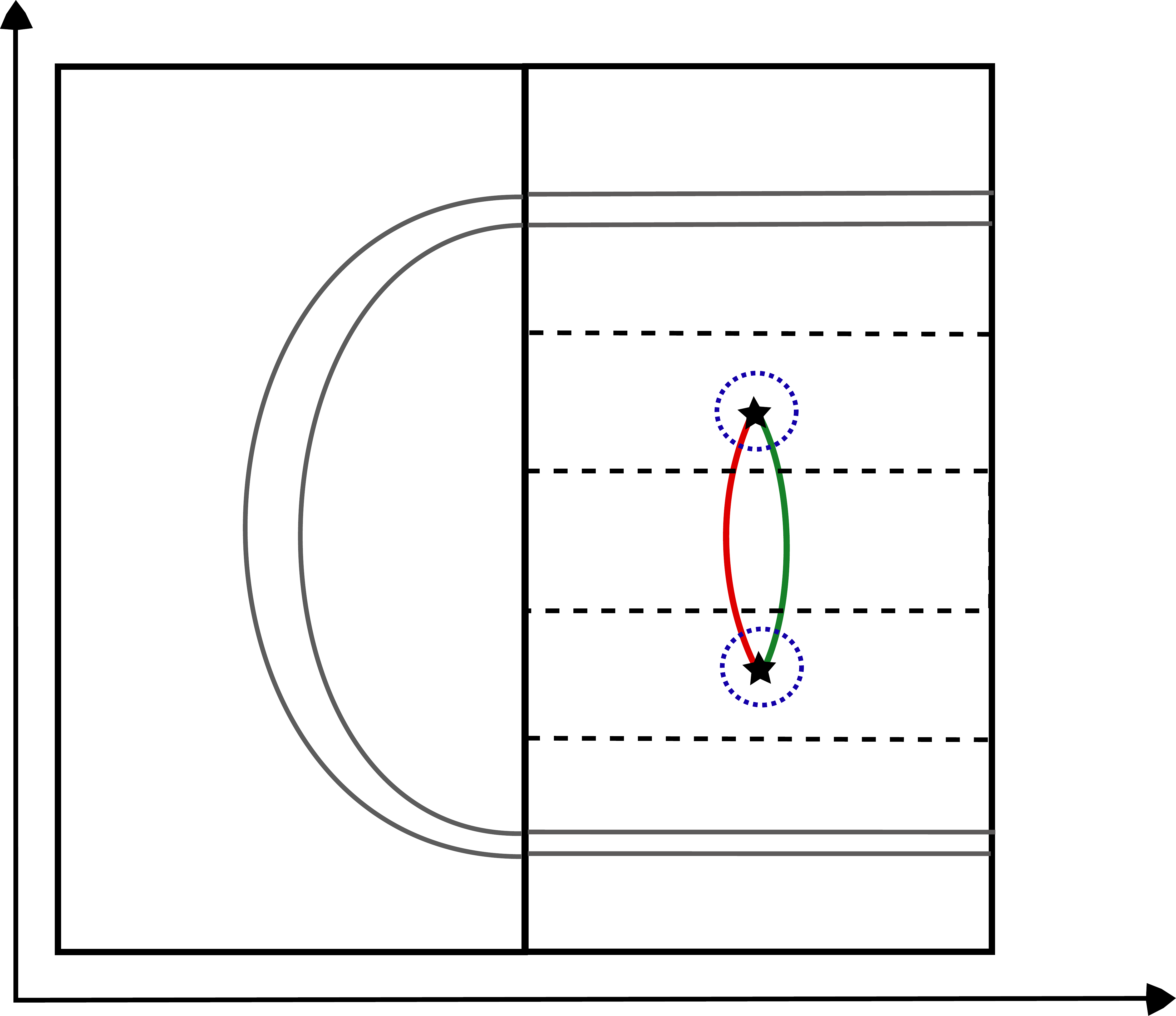}
    \vspace{.3cm}
    \caption{The cusp graphic for the Hatcher-Wagoner family.}
    \label{F: Cusp Graphic}
\end{figure}

In this section we will describe a 2-parameter family interpolating between two Hatcher-Wagoner families. The first family has an empty Cerf graphic and consequentially gives a pseudo-isotopy that is isotopic to an isotopy. The second family generates $[\beta^2] \in \pi_0 \Diff(S^4)$. Although it is already known that $[\beta^2] = [\id]$ \cite{BG1, GH}, it has not yet been observed how the loop of spheres $\mR_t$ that generates $[\beta^2]$ or the Hatcher-Wagoner family interpolates to the trivial family. Here, we provide one such way. Figure \ref{F: Cusp Graphic} provides the \textit{cusp graphic} for the 2-parameter family $(q_{s,t}, \xi_{s,t})$ we build in this section.

\begin{definition}
    For $z\in B^2$, let $(q_z, \xi_z)$ be a generic $2$-parameter family of generalized Morse functions and gradient-like vector fields on $X\times [0,1]$. Let $\Gamma\subset B^2$ denote the set of points such that for $z\in \Gamma$, $q_z$ has a birth-death singularity. The \textit{cusp graphic} is the closure of $\Gamma$. The \textit{handle-slide graphic} is the set of points $z\in B^2$ such that $-\xi_z$ has a flow line from a critical point of index $i$ to one of index $j$ where $i\leq j$. 
\end{definition}

\begin{notation}
    For the pair $(q_{s,t},\xi_{s,t})$, for $(s,t)$ in our parameter space, we will denote the $n$-th critical point of index $i$ $p_n^i(s,t)$. The ascending and descending spheres in an appropriate level set for an index 3 critical point $p^3_n(s,t)$ will be denoted by $\mR^*_n(s,t)$ and $\mR_n(s,t)$ respectively while the ascending and descending spheres for an index 2 critical point $p_m^2(s,t)$ will be denoted by $\mG_m(s,t)$ and $\mG_{m*}(s,t)$ respectively. We will use $\mR(s,t)$, $\mR^*(s,t)$, $\mG(s,t)$, and $\mG_*(s,t)$ to denote the union of all of the ascending/descending spheres in a particular level set.
\end{notation}

\subsection{Loops of spheres.}
 We start with the description of the family that separates region $A$ from region $B$. This is a 2-eyed family; as $t$ increases, from $t=0$ to $t=\frac{1}{8}$ a pair of standard births occurs. From just after $t = \frac{1}{8}$ to $t=\frac{3}{8}$, the family is stationary. From $t=\frac{3}{8}$ to $t = \frac{5}{8}$, the attaching spheres for the 3-handles will dance around. This will be explicitly described below. After the loop of spheres returns, the family is stationary from $t=\frac{5}{8}$ until $t=\frac{7}{8}$. Then the handles cancel, passing through two standard deaths.

 In order to elucidate the evolution of this family, we first need to fully describe the loop of attaching spheres. We first give a construction that specifies the entire loop. Then we will use this to describe the loop of spheres $\mR(1/3, t)$ for $\frac{3}{8}\leq t \leq \frac{5}{8}$. It will be clear from the description that this family of attaching maps can be canceled, yeilding the $A$-block. We will then detail the process of extending this family across the $B$-block. Once we have the description of the family at $s=\frac{2}{3}$, we give the necessary modifications to the family to reduce it to a concatenation of the family in Proposition \ref{P: 12 FW system} with itself.

 \begin{construction}\label{Construction: Spinning}
     Let $\Sigma$ be an embedded compact manifold of genus $g$ with one boundary component embedded in $X^4$. For each $0\leq i \leq 3$ let $\mC_i = \{c_1^i,..., c_g^i\}$ be a collection of disjointly embedded framed disks satisfying the following conditions:
     \begin{enumerate}
        \item For each $i$, the disks $c_k^i\in \mC_i$ only intersect $\Sigma$ on along their boundaries, $c_k^i\cap \Sigma = \partial c_k^i$, and the framing induced by $\Sigma$ on the $D^2$-bundle along $\partial c_k^i$ agrees with the framing of $c_k^i$ restricted to the boundary.
        
        \item Each collection $\mC_i$ can be used to perform embedded surgery on $\Sigma_g$ and the result surgers $\Sigma_g$ down to an embedded disk $\Sigma^{\mC_i}$.
         
        \item For $0\leq i\leq 3$ and $(i+1)\mod 4$, the collections $\mC_i$ and $\mC_{i+1}$ are dual collections in the sense that if $k\neq l$, then $c_k^{i}\cap c_{l}^{i+1}=\varnothing$, and for every $1\leq k\leq g$, $|c^i_k\cap c^{i+1}_k| = |\partial c^i_k\cap \partial c^{i+1}_k| =1$. 
     \end{enumerate}
    Following the terminology of \cite{FQ} (see also \cite{ST}), the collection of disks in $\mC_i$ are called \textit{caps} for the surface. Assumption 3 on the collections $\mC_i$ ensures that for each $i$, and $(i+1)\mod 4$, the pair $\mC_i$ and $\mC_{i+1}$ define an isotopy from $\Sigma^{\mC_i}$ to $\Sigma^{\mC_{i+1}}$ that is supported in a neighborhood of $\Sigma \cup \mC_i\cup \mC_{i+1}$ (see Sect 2.3 of \cite{FQ} or Lemma 2.6 of \cite{ST}). Thus we get an isotopy $\phi_1$ from $\Sigma^{\mC_0}$ to $\Sigma^{\mC_{1}}$, $\phi_2$ from $\Sigma^{\mC_1}$ to $\Sigma^{\mC_{2}}$, $\phi_3$ from $\Sigma^{\mC_2}$ to $\Sigma^{\mC_{3}}$ and $\phi_4$ from $\Sigma^{\mC_3}$ to $\Sigma^{\mC_{0}}$. The concatination of $\phi_1,..., \phi_4$ defines a loop of embeddings $\phi(\Sigma,\mC_0, \mC_1, \mC_2, \mC_3)$ based at $\Sigma^{\mC_0}$ that is supported near $\Sigma\cup \mC_0\cup \mC_1\cup \mC_2\cup \mC_3$ \qed
 \end{construction}

\begin{construction}[\textbf{Constructing $\mR(1/3,t)$ for $3/8 \leq t \leq 5/8$}]\label{Ex: f-12 w-12}
    Let $\#^2S^2\times S^2$, $\mR^{std}= \mR^{std}_1\cup \mR^{std}_2$, and $\mG = \mG_1\cup \mG_2$ be the middle level set configuration for $(q_{1/3, 3/8}, \xi_{1/3,3/8})$ with $\mR$ the descending spheres for the 3-handles and $\mG$ the ascending spheres for the 2-handles. We will describe $\mR(1/3,t)$ for $3/8\leq t \leq 5/8$ by building a loop of embedded spheres based at $\mR^{std}$. 
    
    First, perform a standard finger move along an arc that goes from $\mR_1$ to $\mG_2$. Let $\mR$ be the embedded surface after we performed the finger move. Let $f_{1,2}$ and $S_{1,2}$ denote the finger disk and finger sphere respectively associated with the finger move. Now using the arc $f_{1,2}\cap \mG_2$, we tube $\mR$ to itself. This is done first by thickening the arc $f_{1,2}\cap \mG_2$ with disks of arbitrarily small radius contained in the fibers of the normal bundle of $\mG_2$. As $\mR$ is transverse to $\mG$, the disks at the endpoints are subsets of $\mR$ that contain the intersection points. We can form $\Sigma$ by removing the interiors of the disks in $\mR$ and gluing the boundary of the thickened tube along $f_{1,2}\cap \mG_2$. 
    
    The process of tubing $\mR$ to itself as described above removes the intersections with $\mG_2$ and forms the embedded surface $\Sigma$. This construction also converts the disk $f_{1,2}$ into the cap $c^0$ for $\Sigma$. Moreover, there is an obvious dual cap $c^1$ given by taking one of the small normal disks used to create the tubing. This normal disk will intersect $\mG_2$ in a singular point. See Figure \ref{F: Example_f_12 1} for an illustration.
    
    The finger move creates half of the caps needed to define the loop. The remaining caps for this example will be constructed by modifying $c^0$ and $c^1$ as follows. The cap $c^2$ is constructed from $c^0$ by connecting the disk $c_0$ to the finger sphere $S_{1,2}$ using an arc $\gamma$ in the complement of $\Sigma \cup c_1 \cup \mR_2$ with an endpoint in $c_0$ and the other on $S_{1,2}$. Finally, we build $c_3$ from $c_1$ by tubing $c_1$ a dual sphere of $\mG$ that is disjoint from $\mR_2^{std}$. Figure \ref{F: Example_f_12 2} shows the schematic for these 4 caps.

    Now the isotopy based at $\mR^{std}$ is given concatenating the path $\psi_t$ that goes from $\mR^{std}_1$ to $\Sigma^{c^0}$ with the loop $\phi_t(\Sigma, c^0, c^1, c^2, c^3)$ then with $\psi_{1-t}$. The schematic of the loop for $\mR_t$ is shown in Figure \ref{F: f_12 w_12 loop.}. After an appropriate reparameterization, we have that the index $3$ descending spheres for $(q_{1/3,t}, \xi_{1/3,t})$ in the level set above the two index $2$ critical points are given by $\mR(1/3,t)= \mR_t$. \qed

    \begin{figure}
        \centering
        \includegraphics[width=0.9\linewidth]{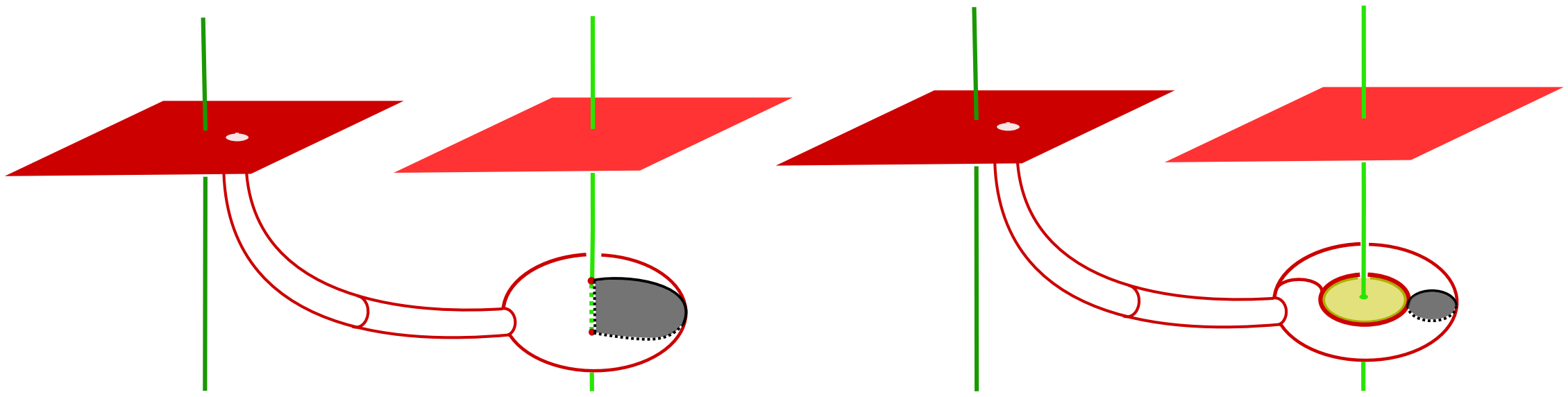}
        \caption{The left hand figure shows the embedding $\mR$ (dark red) which is the result after $\mR_1^{std}$ does a finger move through $\mG_2$ (light green) and the finger disk $f_{1,2}$ (in black). The Figure on the right shows the genus 1 surface obtained from tubing $\mR$ to itself together with the caps $c^0$ (in black) and $c^1$ (in yellow). }
        \label{F: Example_f_12 1}
    \end{figure}

    \begin{figure}
        \centering
        \labellist 
        \small
        \hair 10pt
        \pinlabel $\mG_1$ at 700 1550
        \pinlabel $\mG_2$ at 2250 1550
        \pinlabel $\Sigma$ at 20 850
        \pinlabel $\mR_2$ at 3220 1250
        \pinlabel $S_{1,2}$ at 1700 20
        \pinlabel $\gamma$ at 2000 90
        \endlabellist
        \includegraphics[width=0.7\linewidth]{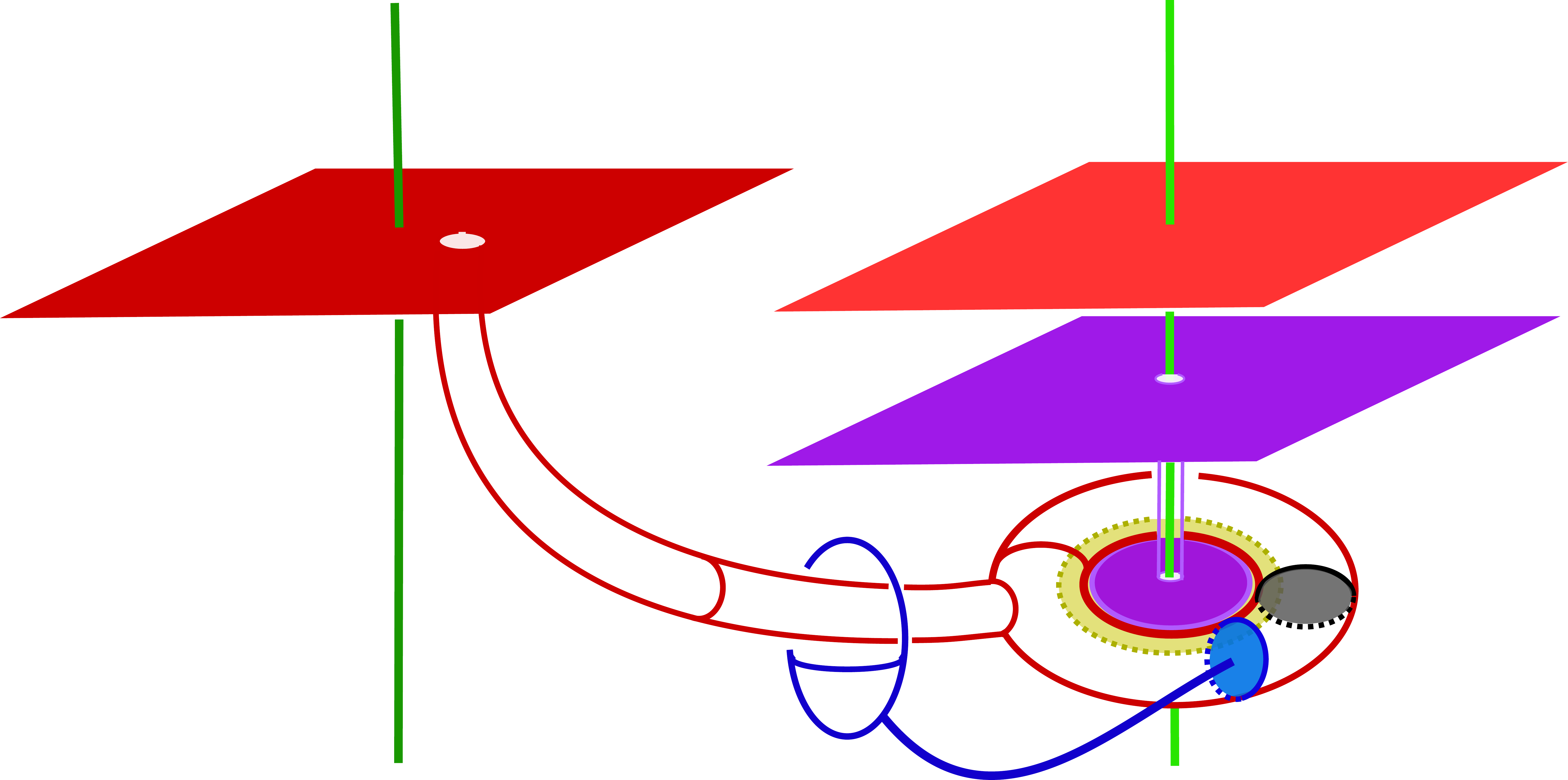}
        \caption{Shown is the surface $\Sigma$ together with the 4 caps $c^0$, $c^1$, $c^2$, $c^3$. Note that $c^3$ (in purple) is built from a parallel copy of $c^1$. Therefore $c^1$ is only partially visible, as it lies below $c^3$ in the figure. The cap $c^2$ (in blue) is built by tubing a parallel copy of $c^0$ (in black) to the finger sphere $S$ along the arc $\gamma$. As the ambient dimension is 4, the arc $\gamma$ does not intersect $\Sigma$.}
        \label{F: Example_f_12 2}
    \end{figure}
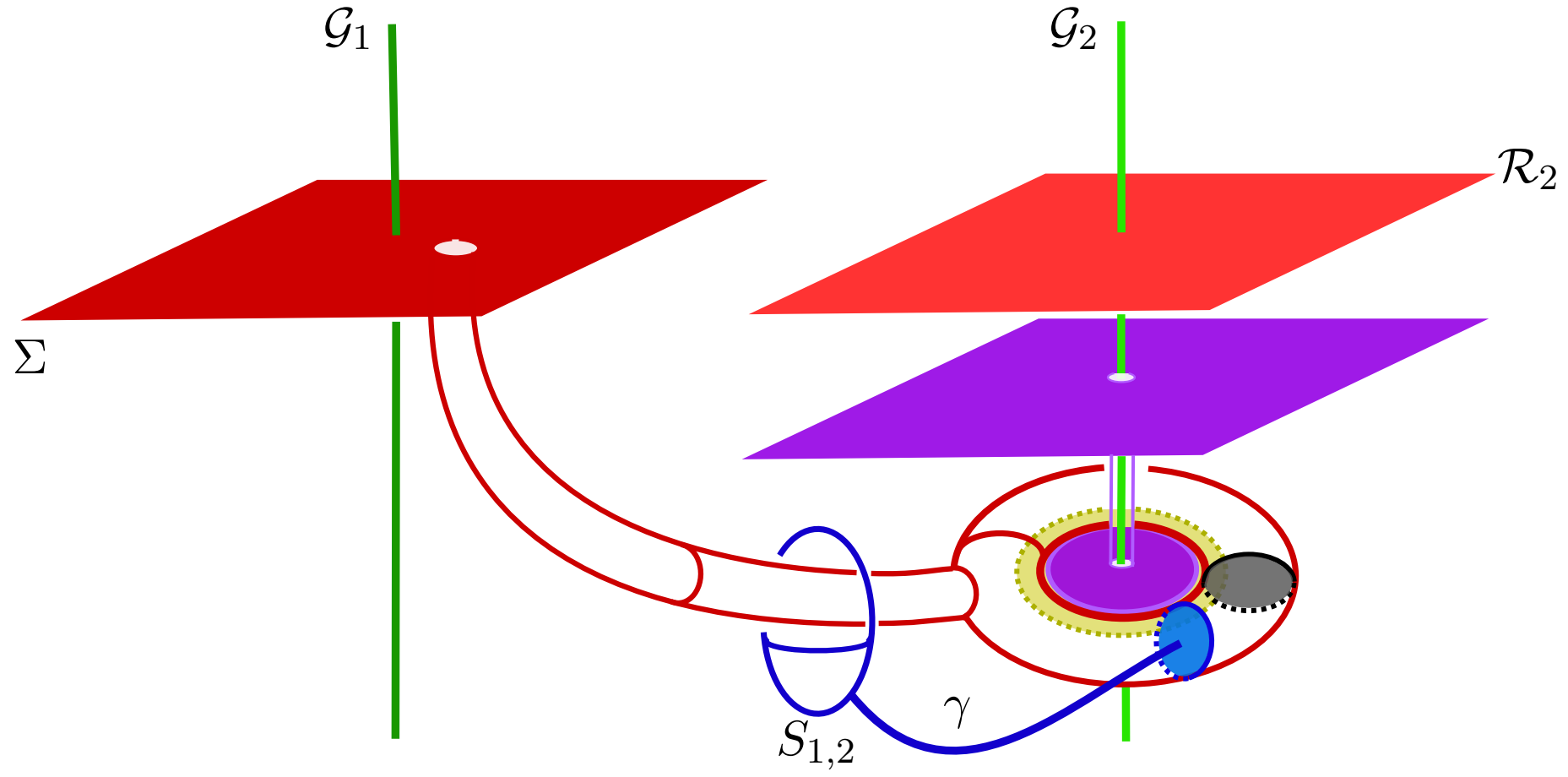

    \begin{figure}
        \centering
        
        \labellist 
        \small
        \hair 10pt
        \pinlabel $\psi_t$ at 3150 3550
        \pinlabel $\psi_{1-t}$ at 3150 1750
        \pinlabel $\phi_{1}$ at 6800 1850
        \pinlabel $\phi_{2}$ at 8200 1850
        \pinlabel $\phi_{3}$ at 8200 2850
        \pinlabel $\phi_{4}$ at 6800 2850
        \endlabellist
        
        \includegraphics[width=\linewidth]{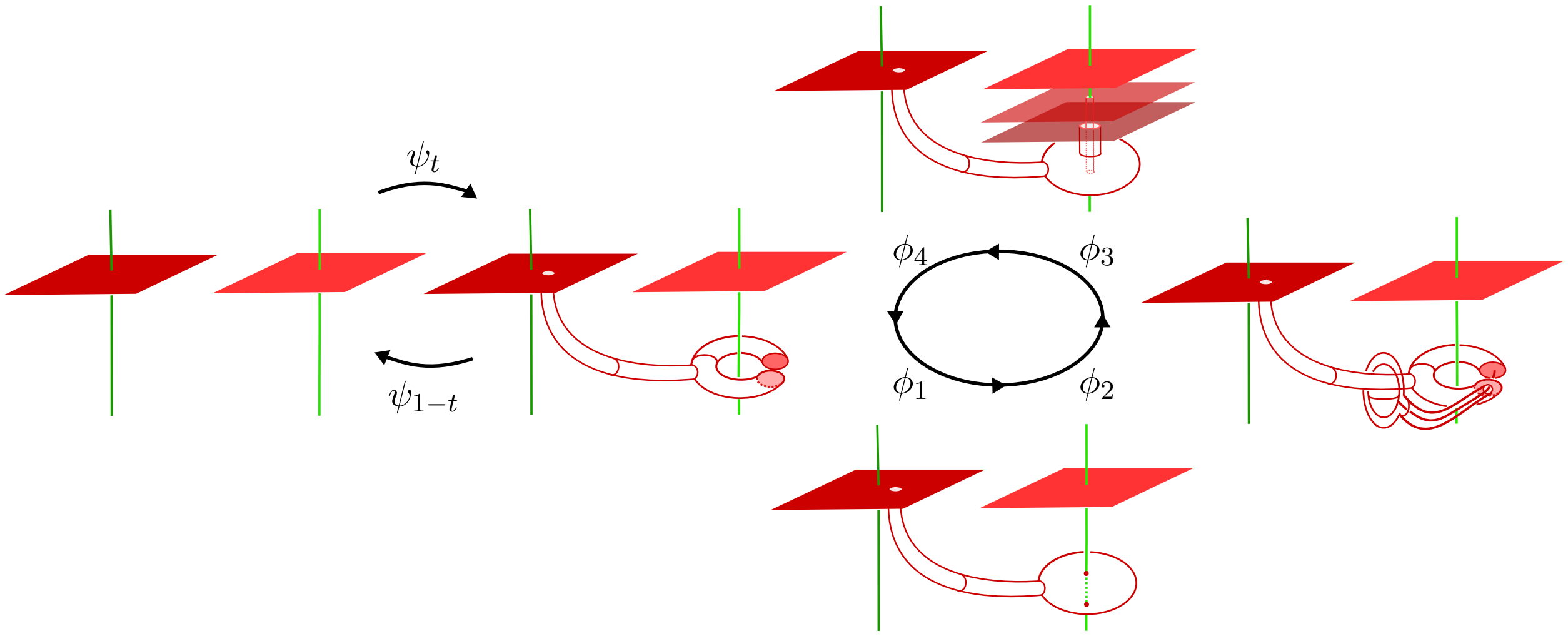}
        \caption{A shown is the loop $\mR_t = (\mR_{1}\cup \mR_2)_t$ based at $\mR^{std}$. The loop is built using Construction \ref{Construction: Spinning} with the caps given in Figure \ref{F: Example_f_12 2}. Note that $\mR_2$ is constant for the duration of the loop.}
        \label{F: f_12 w_12 loop.}
    \end{figure}
\end{construction} 

\begin{remark}[Block A]\label{R: Block A}
The loop of spheres given by Construction \ref{Ex: f-12 w-12} has a finger-Whitney system consisting of a single finger disk and a single Whitney disk. Note that the only intersections that occur happen between $\mR_1$ and $\mG_2$. Therefore, the eye consisting of the pair $\mR_1(t), \mG_1$ is in canceling position for the duration of the loop and can be removed \cite{Ce3}. Since $\mR_2$ does not intersect $\mG_1$ at all, removing the eye does not change the intersection locus between $\mR_2$ and $\mG_2$. Therefore, the second eye may also be removed. Both removals can be performed via a relative homotopy of $(q_{1/3,t}, \xi_{1/3,t})$ and this homotopy is what happens in Block A of Figure \ref{F: Cusp Graphic} as $s$ goes from $s= 1/3$ to $s = 0$. The specifics of this homotopy are not important for our argument; all that matters is that the homotopy exists. That this family has the cusp graphic given in block A is left as an exercise.
\end{remark}

\subsection{Families of handle slides.}
Now we will describe the family $(q_{s,t}, \xi_{s,t})$ for $1/3\leq s \leq 2/3,$ $3/8\leq t \leq 5/8$. This is the region inside the dotted lines of block $B$ shown in Figure \ref{F: Cusp Graphic}. Looking at the figure, one sees that there are two lines that cross the region vertically. These lines indicate the parameter values at which a handle slide occurs. The parameter values on the first line, shown in red, correspond to a $\mR_2$ handle slide over $\mR_1$ while the second line corresponds to a $\mG_1$ slide over $\mG_2$. 

Recall that a handle slide is a path of gradient-like vector fields that contains an $i/i$ trajectory for some pair of critical points with index $i$. So our task in describing the family of functions and gradient-like vector fields over the specified region is to deform $\xi_{1/3,t}$ so as to introduce a family of $3/3$ trajectories and a family of $2/2$ trajectories. This means that for every $3/8\leq t\leq 5/8$, we must give an isotopy of $\mR_2(1/3,t)$ in a level set above $p^3_1(1/3,t)$ so that during the isotopy, the circle $\mR^*_1(1/3,t)$ and the sphere $\mR_2(1/3,t)$ have a single transverse intersection.

To specify the handle slide, we first upgrade $\mR_1(1/3,t)$ to a family of framed embeddings $\nu\mR_1(1/3,t)$ that is compatible with $\xi_{1/3,t}$. By framed embeddings we mean an embedding of $S^2\times B^2$ so that $\nu\mR_1(1/3,t)$ is a family of embeddings of $S^2\times B^2$ such that when restricted to $S^2\times \{0\}$ gives $\mR_1(1/3,t)$. This can be done by first specifying a 1-parameter family of Morse coordinates adapted to $\xi_{1/3,t}$ for $p_1(1/3,t)$ \cite[Cor 3.3]{Laudenbach}. Then, working within these coordinates, we can choose the characteristic embeddings; embeddings of $\phi: S^2\times B^2\rightarrow q^{-1}_{1/3,t}(c_1-\epsilon)\cap V(p_1)$, where $V(p_1)$ is the open set containing the Morse coordinates for $p_1(1/3,t)$ and $\epsilon$ is sufficiently small (see \cite{milnorhcob} for details). Finally, we can flow the image of $\phi$ down into $q^{-1}_{1/3,t}(1/2)$. We denote the image by $\nu\mR_1(1/3,t)$.

Now choose a pushoff $\phi(S^2\times \{z\})$ for $0<|z|<\epsilon$ and flow it down. We denote the image of this sphere by $\mR_1'(1/3,t)\subset \nu\mR_1(1/3,t)$. We now choose a thickened arc $\gamma_{1,2} \simeq B^1\times B^2$ with end points $\{-1\}\times B^2 \subset \mR_2(1/3,t)$ and $\{1\}\times B^2 \subset \mR'_1(1/3,t)$. By shrinking $\epsilon$ if necessary, we may assume that the support of $\mR'_1(1/3,t)$ is not all of the sphere and stays arbitartily close to $\mR_1(1/3,t)$ for all $t$. In this way, we can assume that $\gamma_{1,2}$ is a fixed thickened tube for all $t$ and is contained in the complement of both supports. Note that since both $\mR'_1(1/3,t)$ and $\gamma_{1,2}$ are contained in the complement of $\mR_1(1/3,t)$ for all $t$, both can be flowed into a level set above $p_1(1/3,t)$. Here $\mR_1'(1/3,t)$ bounds a 3--ball that intersects $\mR_1^*(1/3,t)$ in a single point. Let $\Gamma_{1,2}$ denote the union of $\gamma_{1,2}$ and the 3--ball above $p_1$, which is again a 3-ball. The boundary of $\Gamma_{1,2}$ is naturally divided into two disks, one contained in $\mR_2(1/3,t)$ and the other given by $$\mR'_1(1/3,t)\setminus (\gamma_{1,2}(\{1\}\times B^2)) \cup \gamma_{1,2}(B^1\times S^1).$$

We can push the disk in $\mR_2(1/3,t)$ across the $\Gamma_{1,2}$ to the other boundary component. This specifies an isotopy of $\mR_2(1/3,t)$ above $p_1$ which is used to deform $\xi_{s,t}$. Since $\Gamma_{1,2}$ intersects $\mR_1^*$, the isotopy does not extend below $p_1$ but creates a $3/3$ trajectory at the intersection point. Below $p_1$ we can see the result of this isotopy is to perform the connected sum of $\mR_{2}(1/3,t)$ and $\mR'_1(1/3,t)$ along $\gamma_{1,2}$, and this is done for all $t$. Hence for $s = 1/2$ we have $\mR_1(1/2,t) = \mR_1 (1/3, t)$ and $\mR_2(1/2, t) = \mR_2(1/3, t) \#_{\gamma_{1,2}}\mR_1'(1/3,t)$.  
\begin{figure}
    \centering
    \labellist
    \small\hair 2pt
    \pinlabel $s=1/3$ at 1600 100
    \pinlabel $s=1/2$ at 4900 100
    \pinlabel $s=2/3$ at 8100 100
    \endlabellist
    \includegraphics[width=\linewidth]{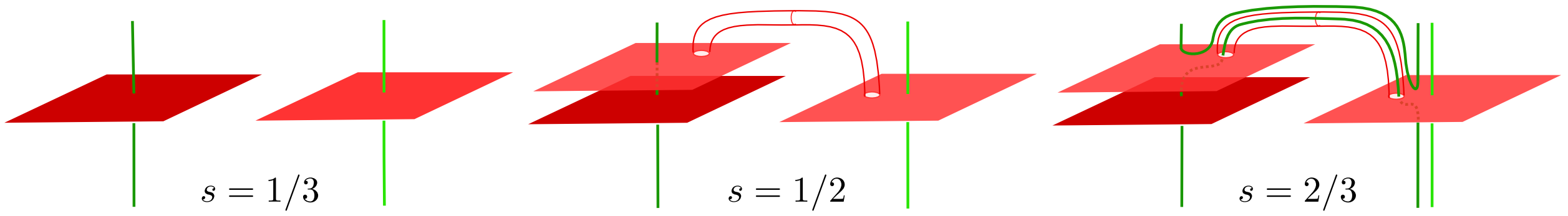}
    \caption{\label{F: exchangehandleslides} Shown here are the configurations of $\mR$ and $\mG$ when $t = 3/8, 5/8$. The left most figure corresponds to $s = 1/3$, the middle to $s  = 1/2$ and the right most figure to $s = 2/3$.}
    
\end{figure}
For $1/2\leq s \leq 2/3$ we repeat this process for $\mG_1$ and $\mG_2$, noting that from the middle level perspective, a $\mG_2$ sliding over $\mG_1$ as index 2-critical points appears as $\mG_1$ sliding over $\mG_2$ as index 3 critical points. Moreover,  we may choose the arc $\gamma^\mG_{2,1}$ to be contained in $\mR_2(1/2,t)$ running from the unique intersection $\mR_2(1/2,t)\cap \mG_1$ to a parallel of $\mG_2$. Also note that in order to perform this handle slide, $p^2_2(s,t)$ must be above $p^2_1(s,t)$. Since this does not affect $\mR$ or $\mG$, as they are above $p^2_2(s,t)$ we can arrange for this to be the case for all $3/8\leq t \leq 5/8$ and $1/3 \leq s \leq 2/3$. This enables us to perform the required $2/2$ handle slide. The configurations of $\mR$ and $\mG$ are shown in Figure \ref{F: exchangehandleslides} for $t = 3/8, 5/8$ and $s=1/3, 1/2, 2/3$ respectively. 

\subsection{Exchange points.}

We now finish describing Block $B$ by describing the family $(q_{s,t}, \xi_{s,t})$ for $1/8\leq t \leq 3/8$ and $5/8 \leq t \leq 7/8$ and for $1/3\leq s \leq 2/3$. Note that these two regions have small circles that contain a star. Outside of the circles, the family will be in a cancellation position. Inside the circles are the regions that consist of parameter values where $p^2_2(s,t)$ rises and is above $p^3_1(s,t)$, and at the star, there is a single flow line from $p_2^2$ to $p_1^3$. In what follows, we will describe the family $(q_{s,t}, \xi_{s,t})$ for the lower portion of Block B, with the understanding that the upper portion is given by an identical process. 

We start with the half-boundary of the lower region that has previously been specified, namely when $s= 1/3, 2/3$ and $1/8 \leq t \leq 3/8$ and $t = 3/8$ and $1/3 \leq s \leq 2/3$. Since $p_1^3$ and $p_2^2$ have no flow lines between them for these parameter values, we apply the Morse rearrangement lemma and homotope $q_{s,3/8}$ relative to $s = 1/3, 2/3$ so that on a subinterval of $1/3 \leq s \leq 2/3$ containing the handle slides, $p_2^2$ rises above $p_1^3$. This homotopy fills the region for $t_0 \leq t \leq 3/8$ and $1/3 \leq s \leq 2/3$, where at $q_{s, t_0}$ is the family where $p_2^2(s, t_0)$ is above $p_1^3(s, t_0)$ on the sub-interval. In Figure \ref{F: Cusp Graphic}, $(1/2, t_0)$ is the point at the top of the lower circle. 

Now we will work in a level set below $p^2_1(s,t)$ but above $p^3_1(s,t)$ and describe how the ascending and descending manifolds behave inside the circle. This will fill in the vector field $\xi_{s,t}$ over the circle with $q_{s,t}$ remaining constant. Once this is done (and up to $(s,t)$-reparameterization), we will have specified $(q_{s,t}, \xi_{s,t})$ over the region $t_1\leq t \leq 5/8$ and $1/3 \leq s \leq 2/3$, where $(1/2, t_1)$ is the bottom point of the lower circle. Observe that the HW family $(q_{s,t}, \xi_{s,t})$ has the property that $p^3_i$ has $\delta_{ij}$ trajectories to $p^2_j$ for the following parameter values: $s=1/3$ and $0\leq t \leq t_1$, $t=t_1$ and $1/3 \leq s \leq 2/3$, and for $s= 2/3$ and $0\leq t \leq t_1$. Therefore there is a 2-parameter family that cancels the eyes \cite[Lemma 3.6]{Laudenbach}. This corresponds to the horizontal lines of cusps at the bottom of Block B in Figure \ref{F: Cusp Graphic}.  

What remains to describe are the ascending and descending manifolds above $p_1^3(s,t)$ and below $p_2^2(s,t)$ over the circle. There we have ${\mG_2}_*$ and $\mR^*_1(s,t)$ as disjoint embedded circles. Since $p^3_2(s,t)$ has a single flow line down to $p^2_2(s,t)$, the descending manifold intersects this level set as a disk $\mR_2(s,t)$ with ${\mG_2}_*(s,t)$ as it's boundary. Similarly, since $p^3_1(s,t)$ has a single flow line to $p_1^2(s,t)$, the ascending manifold $\mG_2(s,t)$ will appear in this level set as a disk with boundary $\mR^*_1(s,t)$.

\begin{figure}
    \centering
    \labellist
    \hair 2pt
    \pinlabel {\small{$\mR_2$}} at -30 950
    \pinlabel {\small{$\mG_{2*}$}} at 300 1180
    \pinlabel {\small{$\mG_1$}} at 2200 950
    \pinlabel {\small{$\mR^*_1$}} at 1950 1180
    \pinlabel {\small{$\gamma_{1,2}$}} at 850 450
    \endlabellist
    \includegraphics[width=\linewidth]{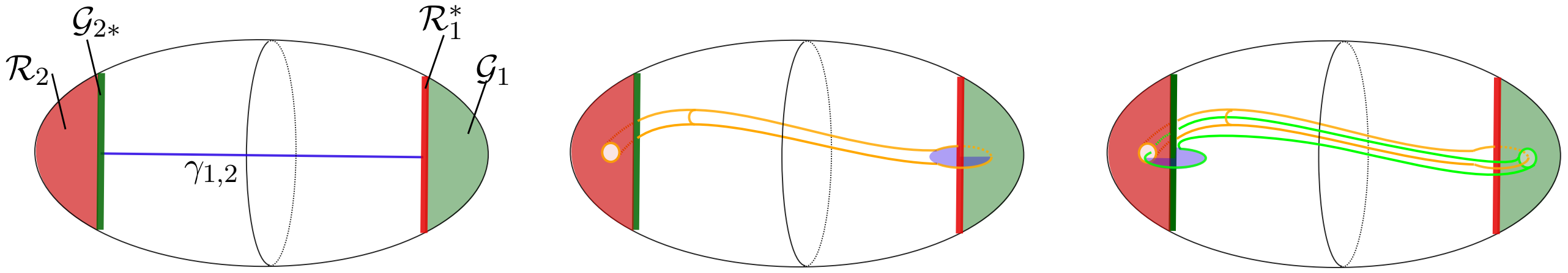}
    \caption{\label{F: 2 above 3 handleslide}This figure shows the $4$-ball neighborhood $B$ of the arc $\gamma_{1,2}$ in a level set above $p^3_1$ and below $p^2_2$. The picture on the left is the initial configurations of $\mR^*_1(s,t)$, $\mR_2(s,t)$, $\mG_1(s,t)$ and $\mG_{2*}(s,t)$ for $t = 3/8$, $s = 1/3$. The middle pictures shows the resulting configurations after the $\mR$-handle slide while the right most picture shows the results of the $\mG$-handle slide. }
\end{figure}

Recall from the previous subsection that the $\mR$-handle slide was constructed above $p^3_1$ and the support for the isotopy is disjoint from $\mG_2$. Similarly, the $\mG$-handle slide was constructed below $p^2_2$ and is disjoint from $\mR_1$. As a result, the isotopies that realize both handle slides are seen in the level set below $p^2_2$ and above $p^3_1$ in which we are working. This is depicted in Figure \ref{F: 2 above 3 handleslide}. Here we see a 4-ball neighborhood $B$ of the arc $\gamma_{1,2}$. The circles $\mR^*_1$ and ${\mG_2}_*$ intersect this neighborhood in a pair of arcs while $\mG_1$ and $\mR_2$ intersect it in ``half''-disks.   

\begin{definition}\label{D: exchange}
    An exchange point $(s_0,t_0)$ is a point in the parameter space such that there exists a single flow line from an index $i$ critical point $p^i$ down to an index (i+1) critical point $p^{i+1}$.
\end{definition}

Working with $\mR^*_1(s,t)$ and ${\mG_2}_*(s,t)$ as in the left of Figure \ref{F: 2 above 3 handleslide}, we can generate an exchange point as follows:
\begin{construction}\label{Con: Exchange Construction}
    Working in the 4-ball neighborhood, bring $\mR^*_1$ and ${\mG_2}_*$ together along $b$ so that they intersect. Now there is an $S^1$ worth of ways to resolve this intersection. This defines a 1-parameter family of configurations for the two arcs where the two arcs ``spin around each other'' (see Definition 4.1 of \cite{BG1}), and a $B^2$-worth of configurations of the two such that on $\partial B^2$ the arcs are spinning and at the center of $B^2$, the arcs are fixed and intersect each other. Moving from the boundary to the center, we consider that the $B^2$ is foliated by cirlces, $B^2 \simeq S^1\times [0,1]/\{(0,z) \sim (0,z')\}$ where for each radial circle $S^1\times \{x\}$, the two arcs have an orbit that moves them closer together as the foliation moves towards the center. As the arcs orbit, they  rigidly move the half disks $\mG_1$ and $\mR_2$ with them. As such, for every $x\in (0,1]$ there is a point $z\in S^1$ where $\mR_2\cap \mR^*_1\neq \varnothing$ and consists of a single point and $\mG_1\cap {\mG_2}_*\neq \varnothing$. Therefore around each $S^1$, there is an $\mR$-handle slide time and an $\mG$-handle slide time.
\end{construction}

\begin{figure}
	\centering
    \labellist
    \hair 2pt
    \pinlabel $1$ at 0 270
    \pinlabel $\nearrow$ [bl] at 0 270
    \pinlabel $\nwarrow$ [br] at 1180 320
    \pinlabel $2$  at 1180 320
    \endlabellist
	\includegraphics[width=.5\linewidth]{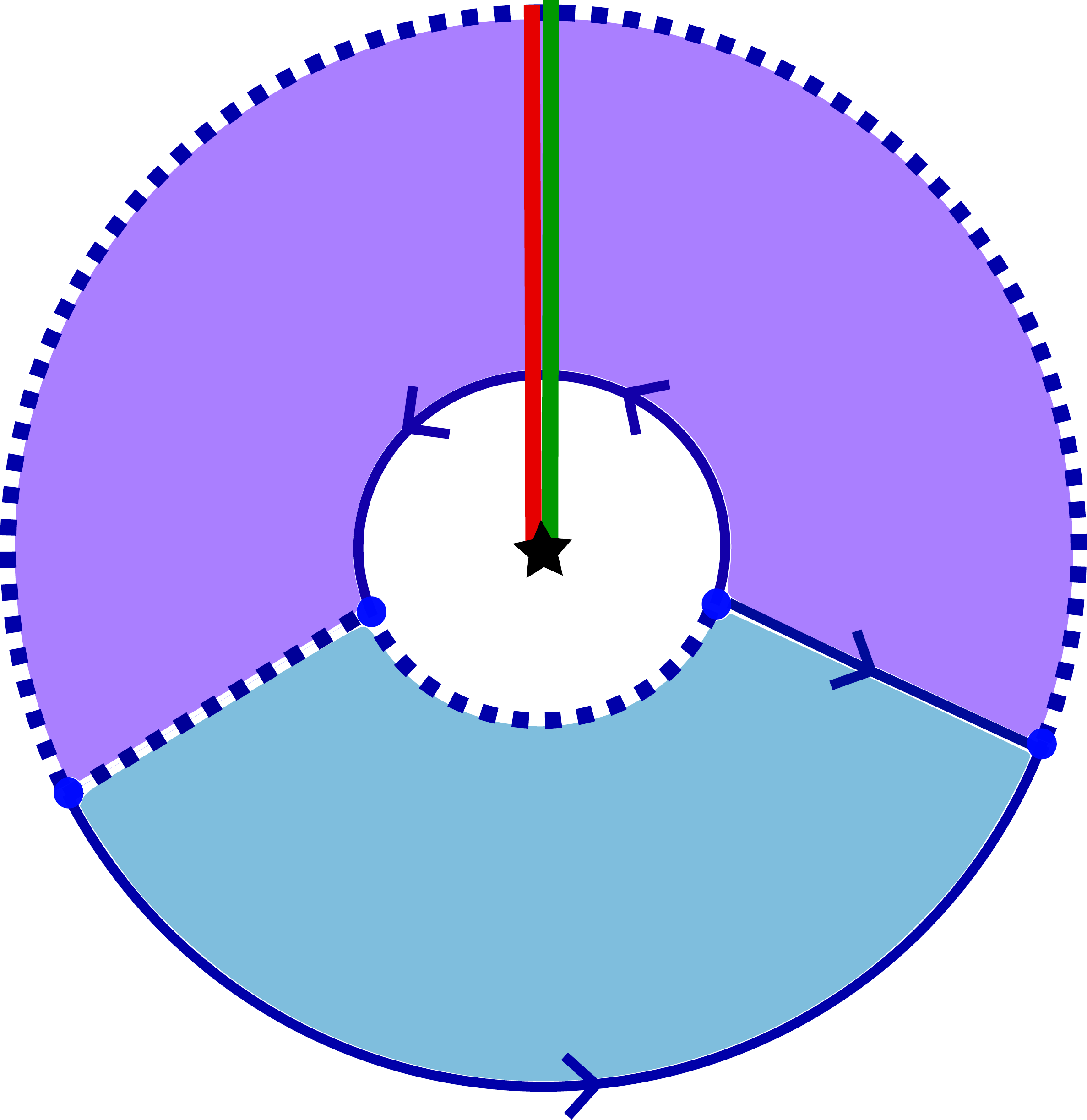}
	\caption{\label{F: def of spinning} This figure shows the annulus of parameters that interpolates from the result of Construction \ref{Con: Exchange Construction} to a configuration of ascending and descending manifolds where $\mR_1^*$ and ${\mG_2}_*$ are stationary from 1 to 2 along the top, and performs the loop of arcs along the bottom. In the interior, $\mG_1$ and $\mR_2$ move rigidly with the arcs, and this motion is extended radially outward using the isotopy extension theorem.}
\end{figure}
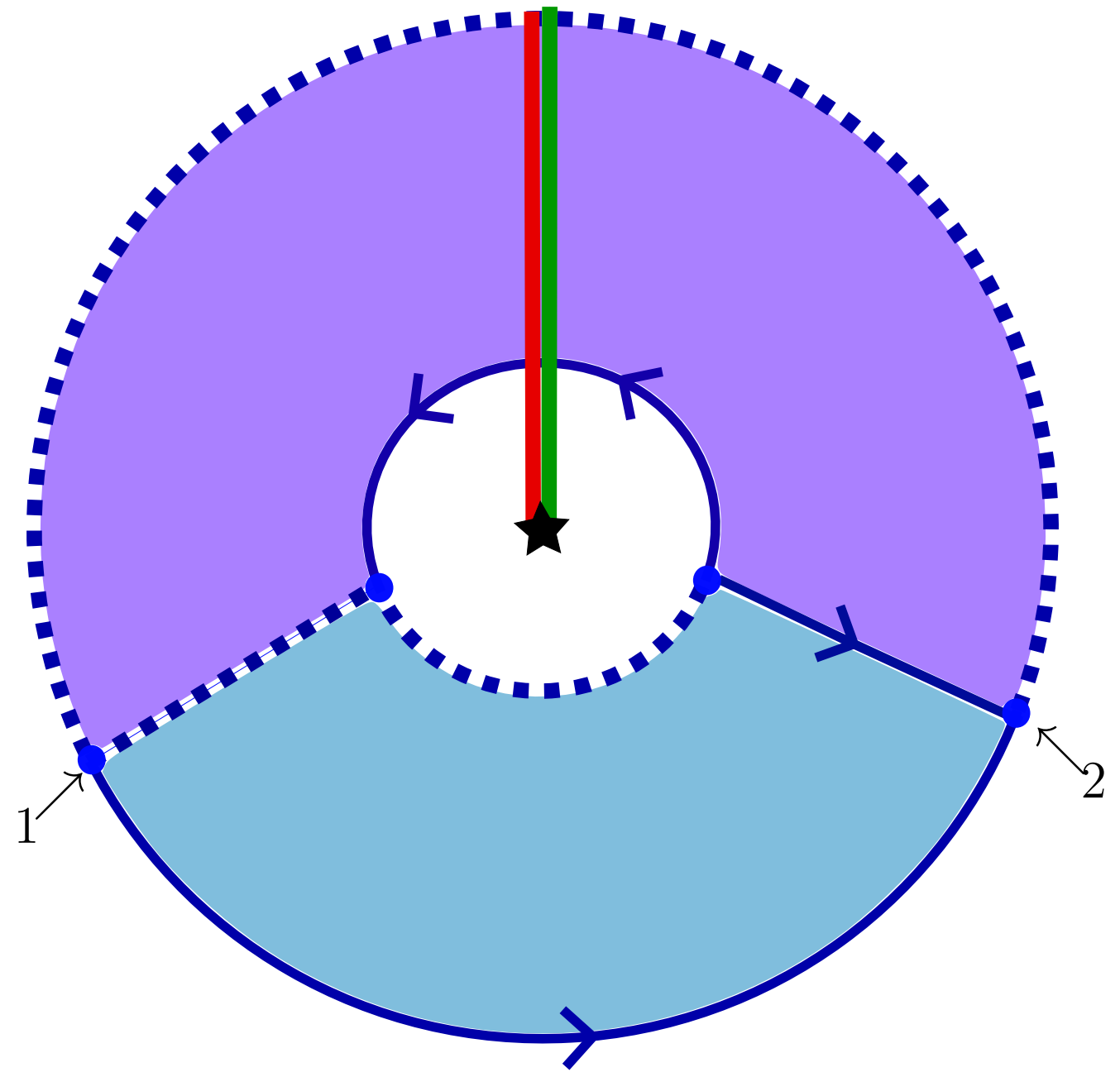

We apply Construction \ref{Con: Exchange Construction} to define $(q_{s,t}, \xi_{s,t})$ over a small $B^2$ subset of the interior of the circular region. We are then left to join this small interior circle with outside boundary that has previously been specified. To do this, we first modify the support of the arc spinning to be disjoint from the handle slides. This is done as follows. 
First suppose that $r(z)$ for $z\in S^1$ parameterizes the loop of embeddings of the two arcs from Construction \ref{Con: Exchange Construction}, and suppose that for $z = e^{i\theta}$, $-\pi/2-\epsilon \leq \theta \leq -\pi/2+\epsilon$, $r(z)$ is constant. We then define a partial homotopy $F_1(z,x)= r({e^{i(\theta -(\theta+\pi/2)x)}})$ for $e^{-\pi/2+\epsilon} \leq z \leq e^{3\pi/2-\epsilon}$ and $0\leq x \leq 1$. We let $F_2(x,z)$ be the null homotopy for $r(z)*r(-z)$, parameterized accordingly to fill in the sector $3\pi/2-\epsilon \leq \theta \leq 3\pi/2 +\epsilon$, $0 \leq x \leq 1$. The domains for these two homotopies are shown in Figure \ref{F: def of spinning}. The dotted lines indicate regions where the arcs are stationary, while the solid lines indicate the support of the loop. The arrows indicate the parametrization direction, so that moving along the line in the given direction performs the loop, while moving opposite to the arrow performs the inverse. Now applying the isotopy extension theorem modifies the disk $\mR_2$ and $\mG_1$ so that the disks perform the handle slides with the arcs stationary. 

\begin{lemma}\label{L: bb exchange pt}
    Let $B$ be the 4-ball neighborhood shown in Figure \ref{F: 2 above 3 handleslide}. Then $B\setminus \mR_1^* \cup {\mG_2}_*$ is diffeomorphism to the model barbell and loop $r(z)$ is the barbell generator. Moreover, applying the ambient isotopy extension theorem to $r(z)$ takes the configuration of $\mG_1\cup \mR_2$ on the left of Figure \ref{F: 2 above 3 handleslide} to the configuration on the right.   
\end{lemma}
\begin{proof}
    See \S 5 of \cite{BG1}.
\end{proof}

It follows from Lemma \ref{L: bb exchange pt} that going from point (1) to point (2) in Figure \ref{F: def of spinning} along the bottom, takes left configuration $\mR_2$ and $\mG_1$ of Figure \ref{F: 2 above 3 handleslide} to the right configuration. This is also the case when moving along the top; however, here we have arranged that $\mR_1^*$ and ${\mG_2}_*$ are fixed. Since this is also the case for our previously defined family, we can match-up the handle slides performed along the top boundary of Figure \ref{F: def of spinning} and those defined in the previous section (see \cite[Chp 2, Lemma 1.2']{HW}). This alows us to extend $\xi_{s,t}$ over the entire circle. Furthermore, we can perform this same construction for the upper exchange point, with all directions reversed. What is left is to then cancel the handles below/above the exchange circles as described above. This completes Block B.

\begin{remark}\label{R: barbell monodromy}
	We point out that the monodromy around the exchange point constructed in Construction \ref{Con: Exchange Construction} is the Barbell diffeomorphism applied to a barbell where one cuff is a parallel $\mR_1$ sphere and the other is a parallel $\mG_2$ sphere. Since the spinning construction of \cite{BG1} is a ``generic unfolding'' of this intersection\footnote{We will not explain this and leave it to the reader to figure out the meaning of this.}, the monodromy around any exchange will be given by the barbell map for some barbell implantation. 
\end{remark}

\subsection{Transferring the finger-Whitney system across the middle.}
\begin{figure}
	\centering
    \labellist
    \hair 2pt
    \pinlabel $\mR_1$ at 150 870
    \pinlabel $\mR_2$ at 1700 870
    \pinlabel $\mG_1$ at 850 1500
    \pinlabel $\mG_2$ at 2400 1500
    \pinlabel $f_{12}$ at 2850 400
    \pinlabel {\Large$\leftarrow$} at 2750 400
    \pinlabel $w_{12}$ at 2600 70
    \pinlabel {\Large$\uparrow$} at 2600 160
    \endlabellist
	\includegraphics[width=\linewidth]{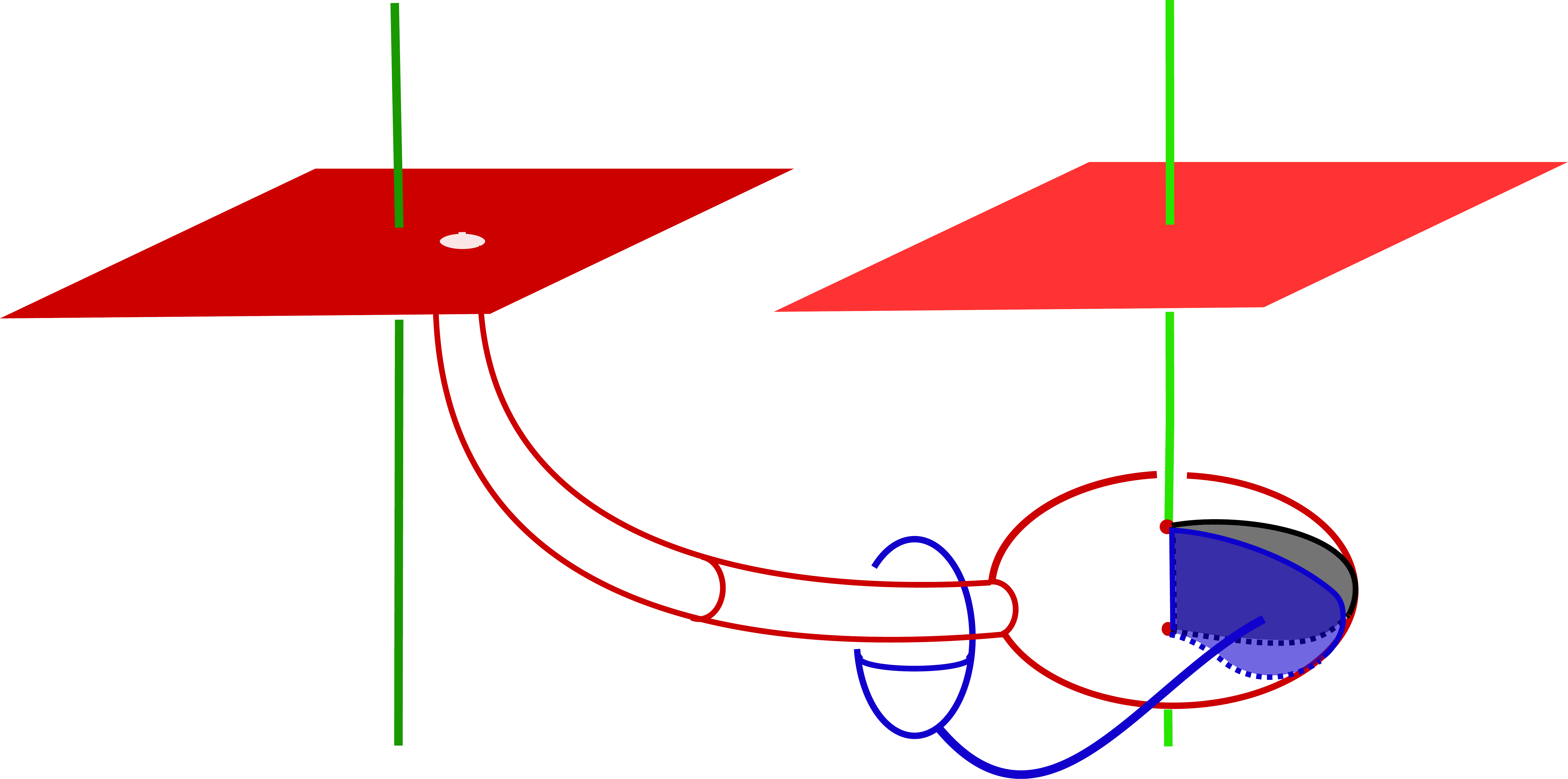}
	\caption{\label{F: FW12} Show is the finger-Whitney system for the loop of spheres $\mR_{1/3,t}$ given by Construction \ref{Ex: f-12 w-12}. }
\end{figure}
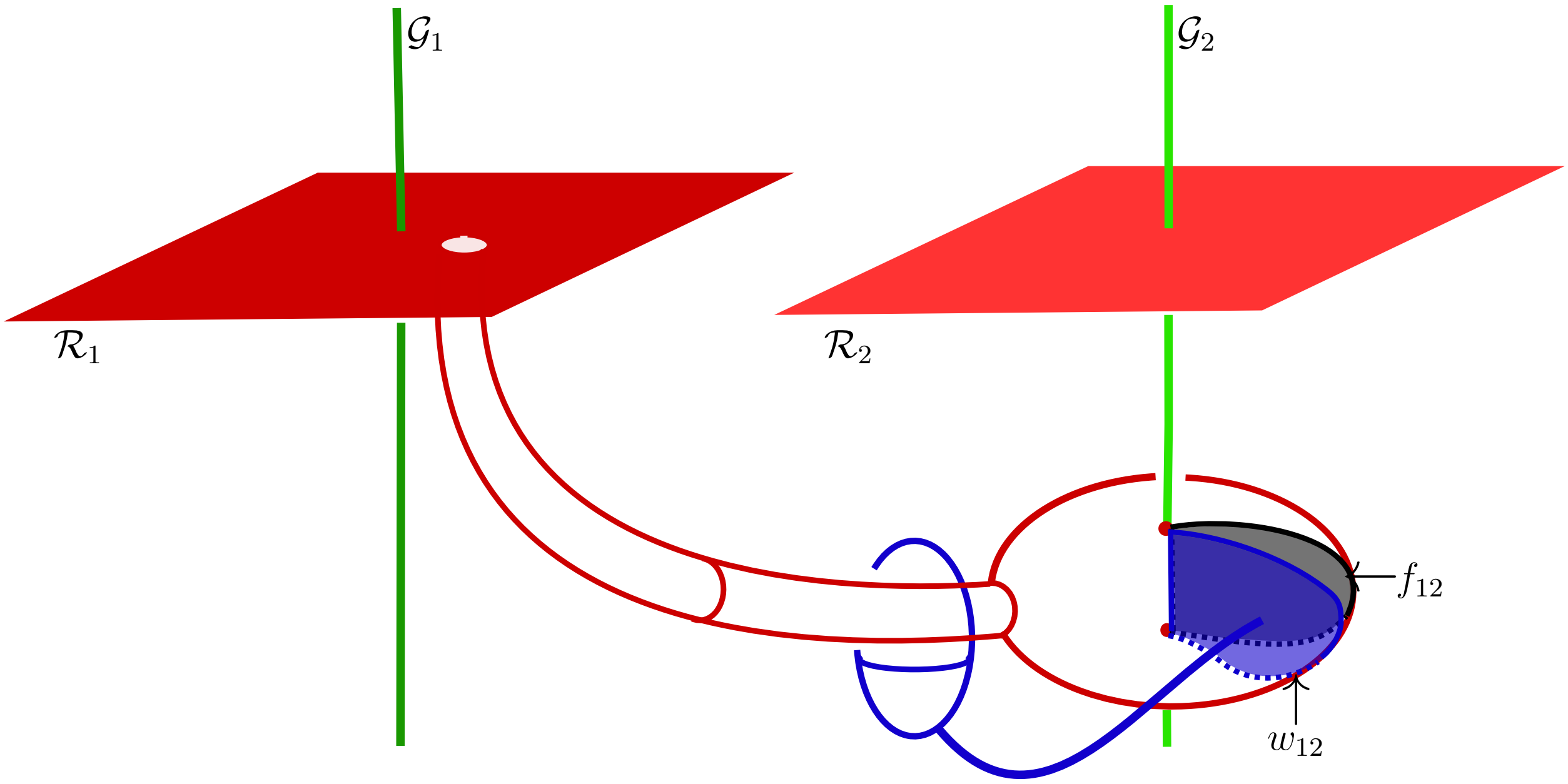

\begin{figure}
	\centering
    \labellist
    \hair 2pt
    \pinlabel $w_{22}$ at 6000 2300
    \pinlabel $w_{21}$ at 5700 2100
    \pinlabel $w_{12}$ at 6000 300
    \pinlabel $w_{11}$ at 5700 100
    \pinlabel $f_{22}$ at 2900 2300
    \pinlabel $f_{21}$ at 2700 2100
    \pinlabel $f_{12}$ at 2900 300
    \pinlabel $f_{11}$ at 2700 100
    \endlabellist
	\includegraphics[width=\linewidth]{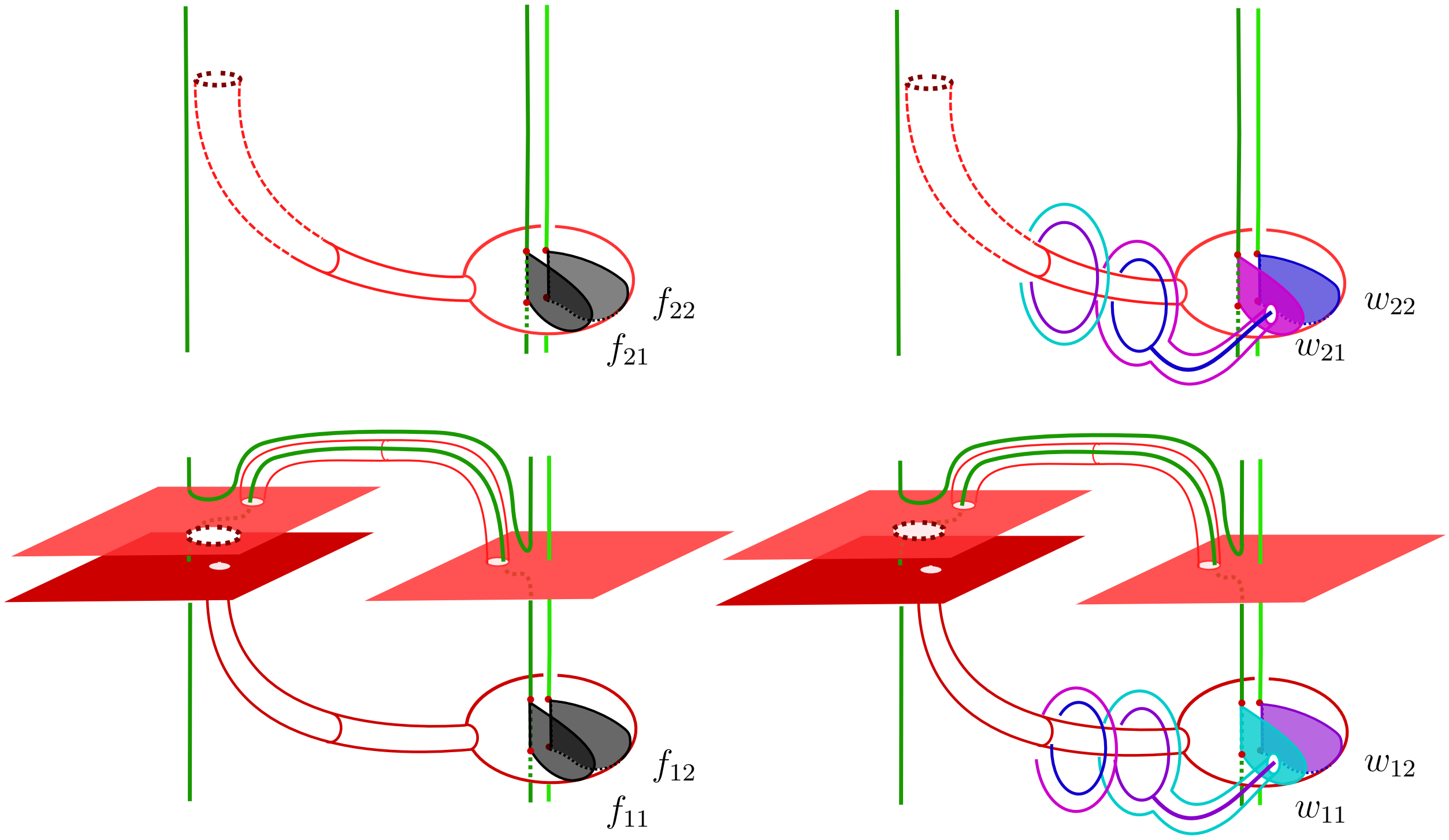}
	\caption{\label{F: FW_12_21_11_22} Shown is the finger-Whitney system after both handle slides have been applied. We have separated the systems so that on the left shows the finger disks and on the right shows the Whitney disks. The dashed lines represent parts of the objects depicted that live out of the 3-dimensional slice shown. For example, the dashed part of the tube in the top two figures connects to the dashed circle in picture below.}
\end{figure}

With Block B specified, we can now analyze the changes that the system goes through crossing the handle slides. As a starting point, we show in Figure \ref{F: FW12} the finger-Whitney system for the family $\mR_{1/3, t}$ from Construction \ref{Ex: f-12 w-12}. Note that the Whitney disk decomposes as the finger disk tubed to the finger sphere. Now when $\mR_2$ slides over $\mR_1$, $\mR_2$ will pick up a parallel copy of every finger disk, and every Whitney disk, since it is being tubed to a ``push off'' of $\mR_1$. Similarly when $\mG_2$ slides over $\mG_1$. Figure \ref{F: FW_12_21_11_22} shows the result of applying the handle slides to the family. The left column of Figure \ref{F: FW_12_21_11_22} shows the finger disks $\mF =\{f_{12}, f_{21}, f_{11}, f_{22}\}$ where the subscripts for $f_{ij}$ indicate the $\mR_i\cap \mG_j$ intersections the disks pair off. The right column shows the corresponding Whitney disk $\mW = \{w_{12}, w_{21}, w_{11}, w_{22}\}$. Here, each column of Figure \ref{F: FW_12_21_11_22} shows two different 3-dimensional slices of $\#^2 S^2\times S^2$. In both of the columns of figures, the dashed lines connect the part of $\mR_2$ that appears in the top frame with the frame below. Since each slide causes the new family to collect parallel copies of the finger and Whitney disks, each Whitney disk for the family at $s=2/3$ decomposes as the corresponding finger disk tubed to a sphere that is the sum of all the finger spheres.

\begin{figure}
\centering
\includegraphics[width= .5\linewidth]{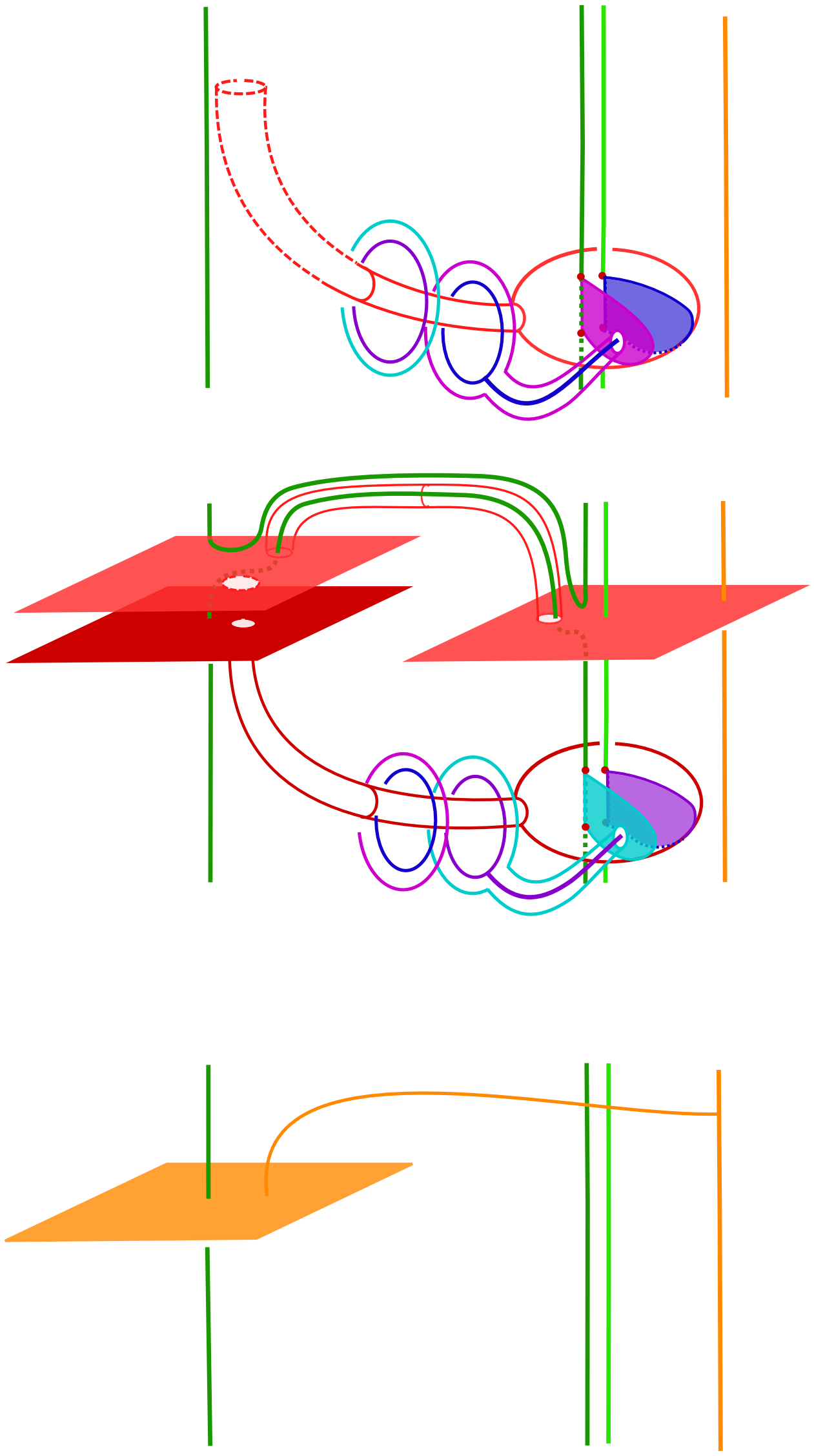}
\caption{\label{F: monodromy} This figure shows the chosen barbell implantation (in orange) for the monodromy. Importantly, this implantation is disjoint from the finger-Whitney system.}
\end{figure}

Note that the spheres $\mR_2$ and $\mG_1$ are not in the correct homology classes to compare the loop with the 1-2 generator. So we must correct for this by applying the monodromy. As discussed in Remark \ref{R: barbell monodromy}, the monodromy for going around an exchange point is given by the barbell diffeomorphism applied to a barbell implantation where one cuff is a parallel of $\mR_1$ and the other cuff is a parallel of $\mG_2$. Let us denote these spheres by $B_1$ and $B_2$ respectively. Since the monodromy is only defined up to isotopy, we are free to choose the parallels so long as $B_1$ is isotopic to $\mR_1$ and $B_2$ is isotopic to $\mG_2$. As such, we choose the spheres $(\{g'\}\times S^2)_1$ and $(S^2\times\{r'\} )_2$ that are disjoint from $(\mF, \mW)$. These cuffs are shown in Figure \ref{F: monodromy} as the orange spheres. 

\begin{definition}
	Let $(\mF, \mW)$ be a finger-Whitney system. We say that $(\mF', \mW')$ is a subsystem if \begin{enumerate}
	\item $\mF'\subset \mF$, $\mW'\subset \mW$,
	\item for any disk $f\in \mF'$, $f \cap \mW\setminus \mW' = \varnothing$ and for any disk $w\in \mW'$, $w\cap \mF\setminus \mF'=\varnothing$. 
	\end{enumerate}
\end{definition}

\begin{lemma}\label{L: FW splitting}
	Let $\mR(t)$ denote a path of embedded spheres in the 4-manifold $X$ such that all the finger moves with the surface $\mG\subset X$ happen before the Whitney moves and let $(\mF, \mW)$ be the finger-Whitney system for $\alpha$. Let $(\mF',\mW')$ be a subsystem of $(\mF, \mW)$. Then the Whitney moves along the disks $\mW'$ commute with the Whitney moves along the disks $\mF \setminus \mF'$. Consequentially, the path $\mR(t)$ factors into a concatenation of two paths with finger-Whitney systems $(\mF', \mW')$ and $(\mF\setminus \mF', \mW\setminus \mW')$.
\end{lemma}
\begin{proof}
	We will work with the track of the isotopy $\tr(\mR(t))\subset X\times [0,1]$. This is the level preserving embedding $(\mR(t), t)$. Note that we can change the heights of the finger moves so that the moves corresponding to $\mF'$ happen before the moves that correspond with the disks $\mF\setminus \mF'$. Now since every disk in $\mW'$ is disjoint from every disk in $\mF\setminus \mF'$, we can extend these disks as in Construction 11.31 of \cite{GGH1} below the times corresponding to the $\mF\setminus \mF'$ finger moves and they will remain Whitney disks. Therefore, we may use them to deform $\tr(\mR(t))$ in a level preserving way by performing the Whitney move along the extended disks until we have removed the intersections created by $\mF'$ finger moves before the $\mF\setminus\mF'$ finger moves.  This defines a homotopy from $\mR(t)$ to a new path $\mR'(t)$ such that $\mR'(t)$ first does the finger moves associated with $\mF'$ then does the Whitney moves associated to $\mW'$ then performs the remaining finger moves then the remaining Whitney moves. The conclusion about the two finger-Whitney systems now follows. 
\end{proof}

\begin{corollary}
	The finger-Whitney system from Figure \ref{F: FW_12_21_11_22} has 4 distinct subsystems in which we may apply Lemma \ref{L: FW splitting}. 
\end{corollary}
\begin{proof}
	As can be see from Figure \ref{F: FW_12_21_11_22}, the subsystems consist of $(\{f_{ij}\}, \{w_{ij}\})$ for $1\leq i,j \leq 2$.  
\end{proof}

Applying Lemma \ref{L: FW splitting}, we split the original system as $(f_{11}, w_{11})*(f_{22}, w_{22})* (f_{12}, w_{12}*(f_{21},w_{21})$. 

\begin{figure}
	\centering
	\includegraphics[width = .7\linewidth]{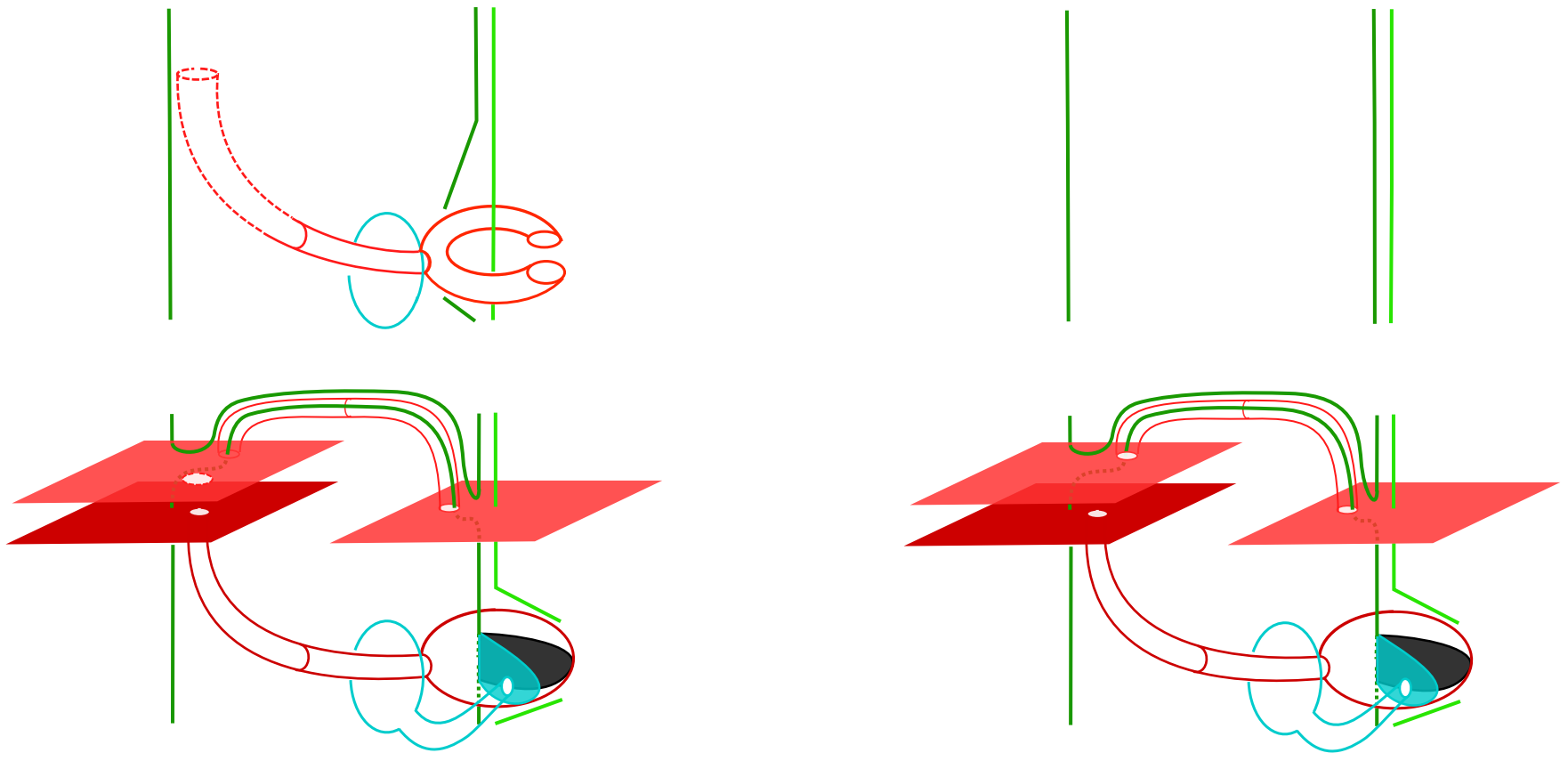}
	\caption{\label{F: Norm11system} The following shows the isotopy normalizing the finger-Whitney system $(\{f_{11}\}, \{w_{11}\})$. }
\end{figure}

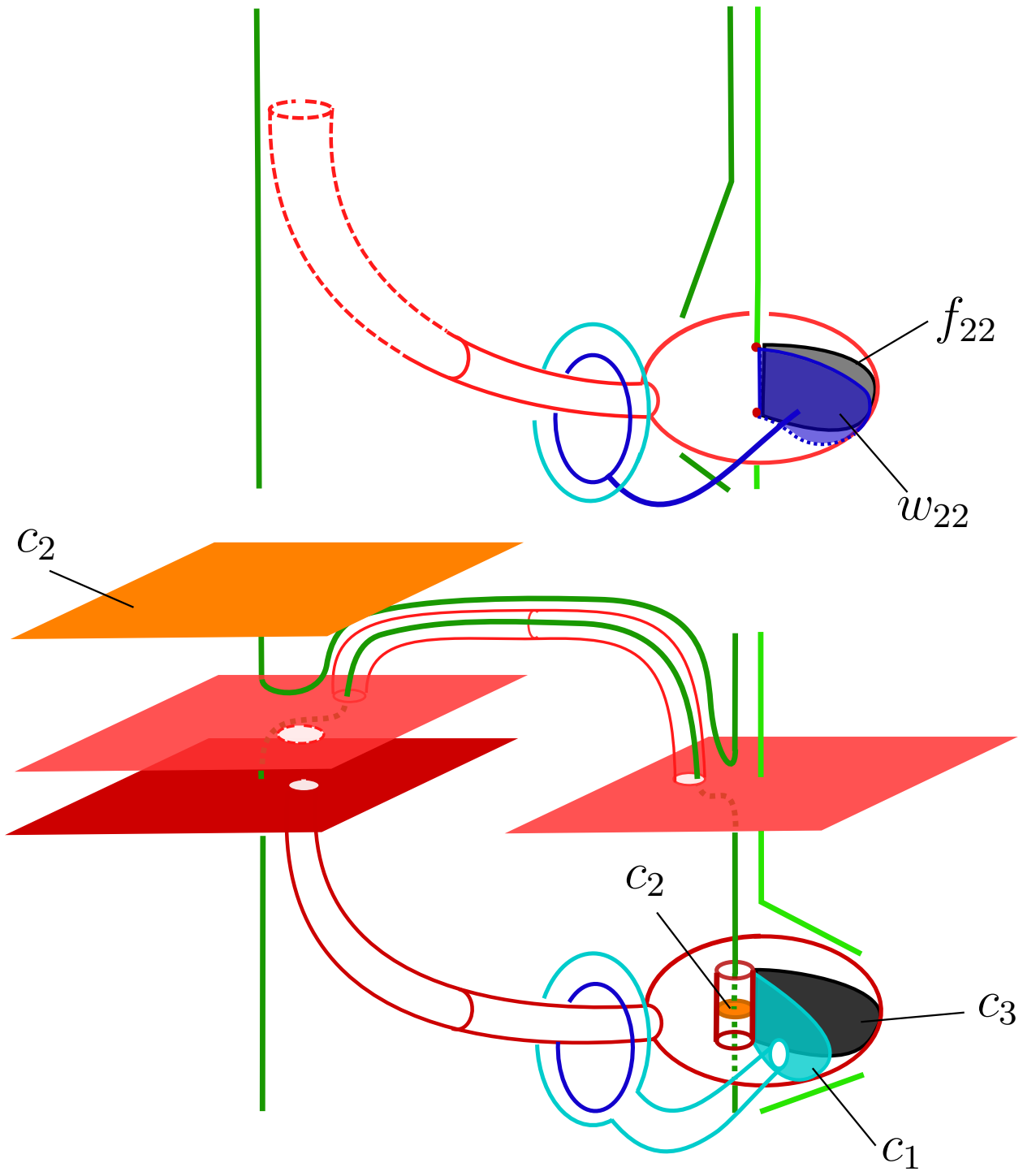
\begin{figure}
	\centering
    \labellist
    \pinlabel $c_1$ at 2800 0 
    \pinlabel {$c_2$} at 100 1900
    \pinlabel {$c_2$} at 2000 850
    \pinlabel $c_3$ at 3100 450
    \pinlabel $f_{22}$ at 3000 2600
    \pinlabel $w_{22}$ at 2900 2000
    \endlabellist
	\includegraphics[width = .5\linewidth]{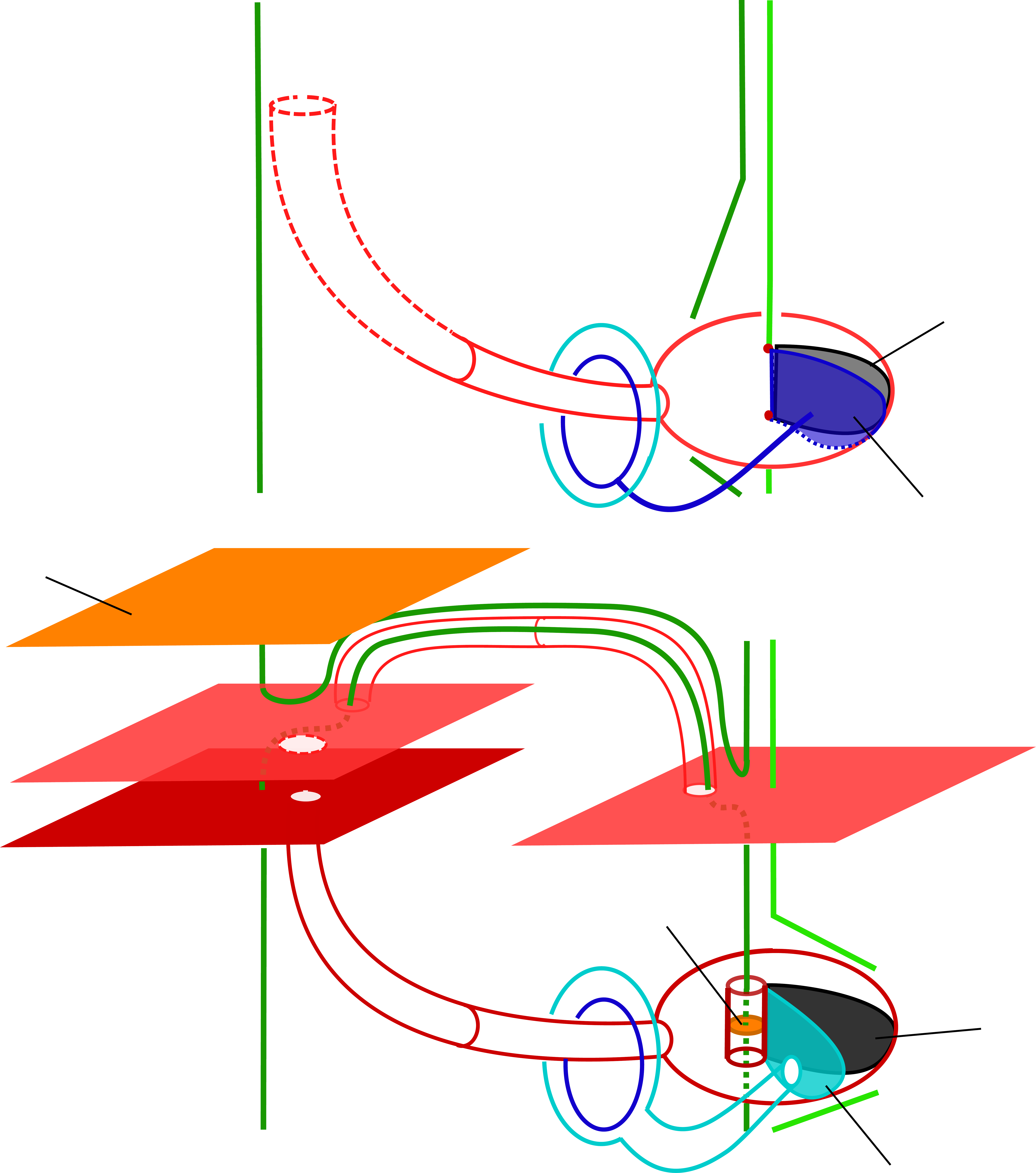}
	\caption{\label{F: LBnormalizingiso} The following shows the 3 caps, $c_1, c_2, c_3$ used to define the first normalizing isotopy for the finger-Whitney system $(\{f_{22}\},\{ w_{22}\})$.}
\end{figure}

\begin{figure}
	\centering
	\includegraphics[width = .7\linewidth]{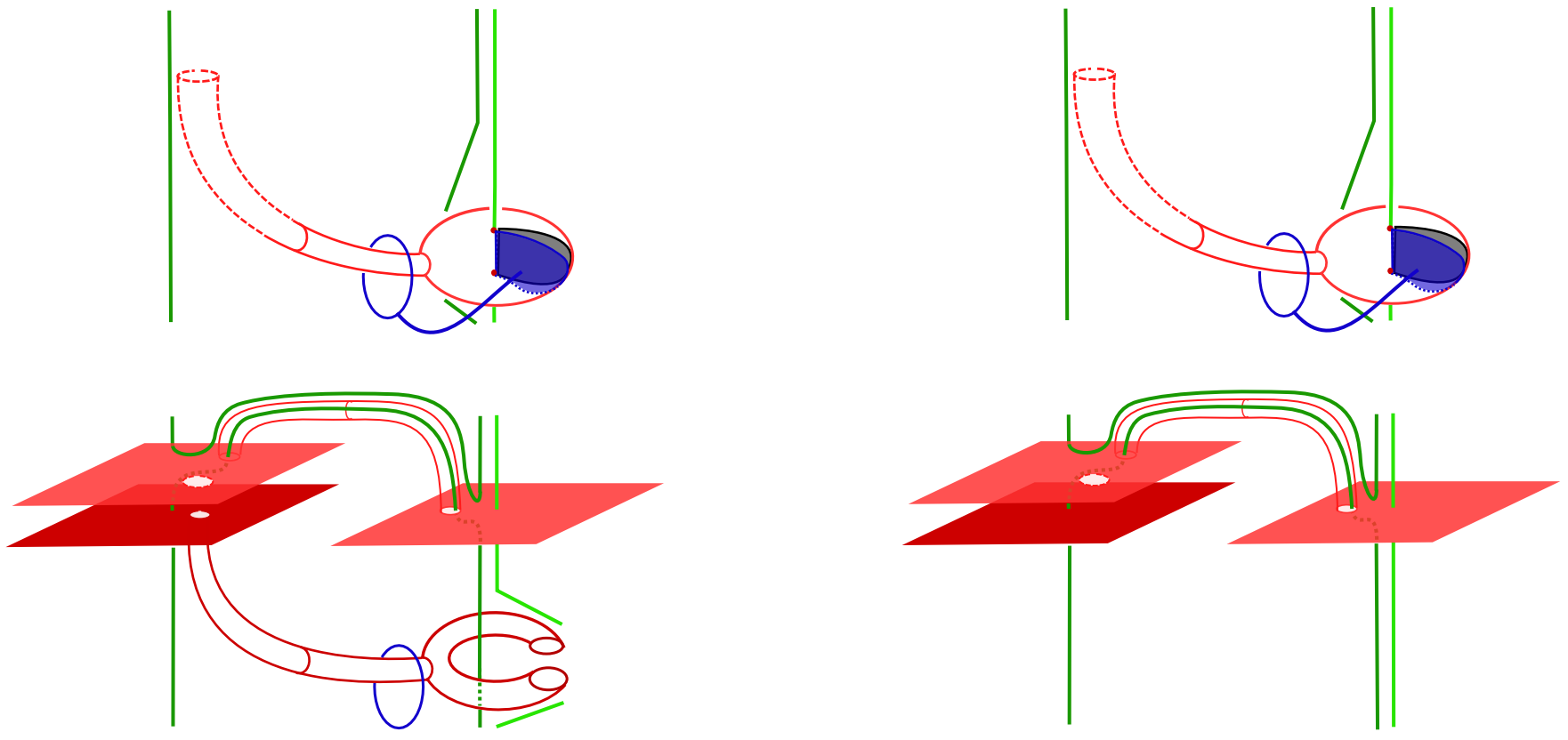}
	\caption{\label{F: Norm22system} The following shows the last normalizing isotopy for $(\{f_{22}\},\{ w_{22}\})$. }
\end{figure}

\begin{lemma}\label{L: Normalization of 11 and 22 systems}
	After rearranging, we can isotope the individual finger-Whitney systems so that $w_{ii}$ decomposes as $f_{ii}$ tubed to the finger-sphere for $f_{ii}$.
\end{lemma}
\begin{proof}
	The proof that $w_{11}$ decomposes as $f_{11}$ tubed to $f_{11}$'s finger sphere is given in Figure \ref{F: Norm11system}. On the left we see the configuration of $\mR_1$ and $\mR_2$ after the $f_{11}$ finger move, but before all other finger moves. The right column of figures shows the obvious isotopy. It is then clear that $w_{11}$ decomposes as claimed. 
	
	To show $w_{22}$ also does, we must apply a larger, non-obvious isotopy. The data for this isotopy is shown in Figure \ref{F: LBnormalizingiso}. Here, we have constructed a genus 1 surface by tubing away the $\mR_1\cap \mG_1$ intersections. If we let $\Sigma$ denote this surface, then surgery on the cap $w_{11}$ results in the configuration of $\mR_1$ just after the Whitney move along $w_{11}$. This is the starting point for the isotopy. We then use the orange cap followed by the $f_{11}$ cap to isotope $\mR_1$ as outlined in Construction \ref{Construction: Spinning} to isotope of $\mR_1$ in the complement of $\mR_2\cup \mG_2\cup \mG_1\cup f_{22}\cup w_{22}$ to surgery on $f_{11}$. Since this isotopy is supported in the complement of $f_{22}\cup w_{22}$, we can  perform an isotopy of $\mR_1$ from just after the Whitney move on $w_{11}$ to after the Whitney move on $w_{22}$ by inserting this isotopy and its inverse. To be more specific, we deform $\mR_1$ in the level sets of $\tr(\alpha')$ after the Whitney move on $w_{11}$ so that $\mR_1$ isotopes to $\Sigma^{f_{11}}$ and is then constant while $\mR_2$ performs the finger move creating $f_{22}$ and then removes these intersections with $w_{22}$. After this, the inverse isotopy returns $\mR_1$ to the configuration of surgery on $\Sigma^{w_{11}}$.
	
	Shown in the left column of Figure \ref{F: Norm22system} is the finger-Whitney system $(\{f_{22}\}, \{w_{22}\})$ after the applying the constructed deformation. The right column shows the obvious standardization. Clearly $w_{22}$ decomposes as $f_{22}$ tubed to its finger-sphere.
\end{proof}

\begin{lemma}\label{L: monodromy preserves fw}
	For $1\leq i,j\leq 2$, $(\{f_{ij}\}, \{w_{ij}\})$ is disjoint from the support of the monodromy.
\end{lemma}
\begin{proof}
	As each subsystem was disjoint by our choice of barbell implantation, the deformation given in Lemma \ref{L: FW splitting} preserves this. 
\end{proof}

\begin{remark}
	The normalizing isotopy used in Lemma \ref{L: Normalization of 11 and 22 systems} is affected by the monodromy. If we let $\mR(t)$ denote the path of spheres and $\mR'(t)$ the path after the deformation, then the barbell map takes $\mR'$ and maps it to a loop with a different base point. This changes the parts of the loop for $\mR'(t)$ that intersect the barbell implantation. However, the finger-Whitney disks and the local isotopies they support remain unchanged as the disks are disjoint from the barbell implantation. 
\end{remark}

Now we can apply the monodromy map to $\mR(t)$ and compare the resulting finger-Whitney systems with the finger-Whitney system for the loop of spheres in Section 3.

\begin{proposition}
	Let $\phi_B$ denote the monodromy diffeomorphism going around either exchange point in block B. Let $\mR'(t)$ denote the modified path that consists of the concatenation of the four finger-Whitney systems $\{f_{11}, w_{11}\}*\{f_{22}, w_{22}\}*\{f_{12}, w_{12}\}*\{f_{21}, w_{21}\}$. Then $\delta([\phi_B(\mR'(t))]) = [\beta^2]$. 
\end{proposition} 

\begin{proof}
	First, we note that between each finger-Whitney system, we have that $\mR_i\cap \mG_j= \delta_{ij}$, and hence the gradient-like vector field is in a canceling position. Thus, the Hatcher-Wagoner family splits into 4 distinct eyes. Over the intervals described by $(\{f_{12},w_{12}\})$ and $(\{f_{21},w_{21}\})$, these eyes can be canceled just as was done over Block A. Over the segment described by $(\{f_{11}\}, \{w_{11}\})$, $\mR_2\cup \mG_2$ are always in a canceling position; thus, this eye can be removed, leaving only $\mR_1\cup \mG_1$. Similarly, we can cancel $\mR_1\cup \mG_1$ over the interval described by $(\{f_{22}\}, \{w_{22}\})$. Thus, we have reduced it to two 1 eyed families. 

Now, by Lemmas \ref{L: Normalization of 11 and 22 systems} and \ref{L: monodromy preserves fw}, we know that $w_{ii}$ decomposes as $f_{ii}$ tubed to the $f_{ii}$ finger sphere. Up to isotopy, the $f_{ii}$ and the finger sphere are contained in a neighborhood of the finger arc from $\mR_i$ to $\mG_i$. As such, we can apply an isotopy to standardize $\mR_i$. Since $S^2\times S^2 \setminus (\mR_i^{std}\cup \mG_i)$ is simply connected, the finger move can be made standard. Hence, the finger-Whitney system  $(\{f_{ii}\}, \{w_{ii}\})$ is equivalent to the finger-Whitney system of Proposition \ref{P: 12 FW system}. Therefore, both systems give rise to loops of spheres $\widetilde{\mR}_i(t)$ such that $\delta([\widetilde{\mR}_i]) = [\beta]$. Merging the eye results in concatenating the two loops of spheres. Since $\delta$ is a group homomorphism, we have that $\delta([\widetilde{\mR}_1]*[\widetilde{\mR}_2]) = 2[\beta] = [\beta^2]$. Since the above modifications are performed (implicitly) with $2$-parameter families of Morse functions and gradient like vector fields, starting from the family generated by$\mR'(t)$ to the family generated by $\widetilde{\mR}_1*\widetilde{\mR}_2$, it follows that $\delta([\phi_B(\mR'(t)]) = [\beta^2]$. 
\end{proof}

\begin{remark}
    Recall that $\pi_0\Diff(X)$ inherits a group structure from composition: $[f]*[g]:=[f\circ g]$, where $*$ is the group operation on $\pi_0\Diff(X)$. Note that function composition is not commutative. Hence $\Diff(X)$ is not an abelian group. For $S^4$, it is known that $\pi_0\Diff(S^4)$ is an abelian group. As such, we use the additive notation $``+''$ to represent the group operation on $\pi_0\Diff(S^4)$. So note that $[\beta^2] = 2[\beta]$, where the superscript refers to function composition, while the $2$ outside refers to the group structure on isotopy classes.
\end{remark}

\begin{corollary}
	$[\beta^2]=[\id_{S^4}]$ 
\end{corollary}
\begin{proof}
 Block A and Block B give a 2-parameter family of Morse functions and gradient-like vector fields interpolating from the Hatcher-Wagoner family with $\phi_B(\mR'(t))$ to the trivial Hatcher-Wagoner family. Hence $\beta^2$ is isotopic to the identity.
\end{proof}

\section{The Double Cross Family}
\begin{figure}
    \centering
    \includegraphics[width=.7\linewidth]{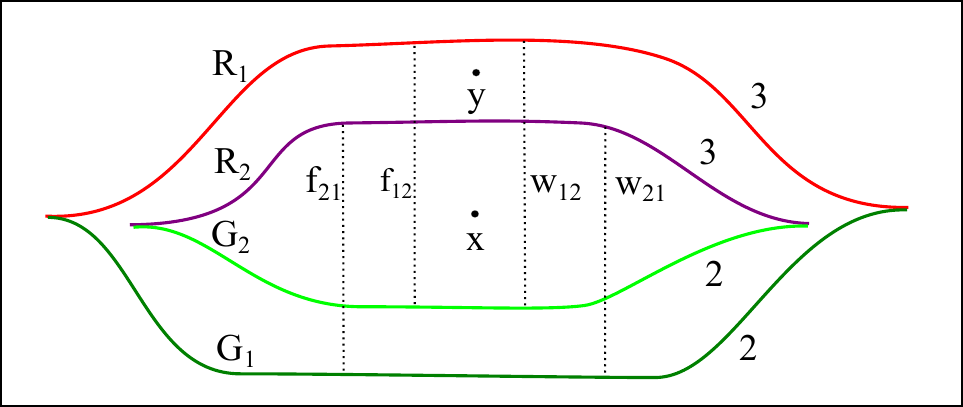}
    \caption{Graphic of the double cross family}
    \label{fig:graphic}
  \end{figure}

\begin{figure}
    \centering
    \includegraphics[width=1.\linewidth]{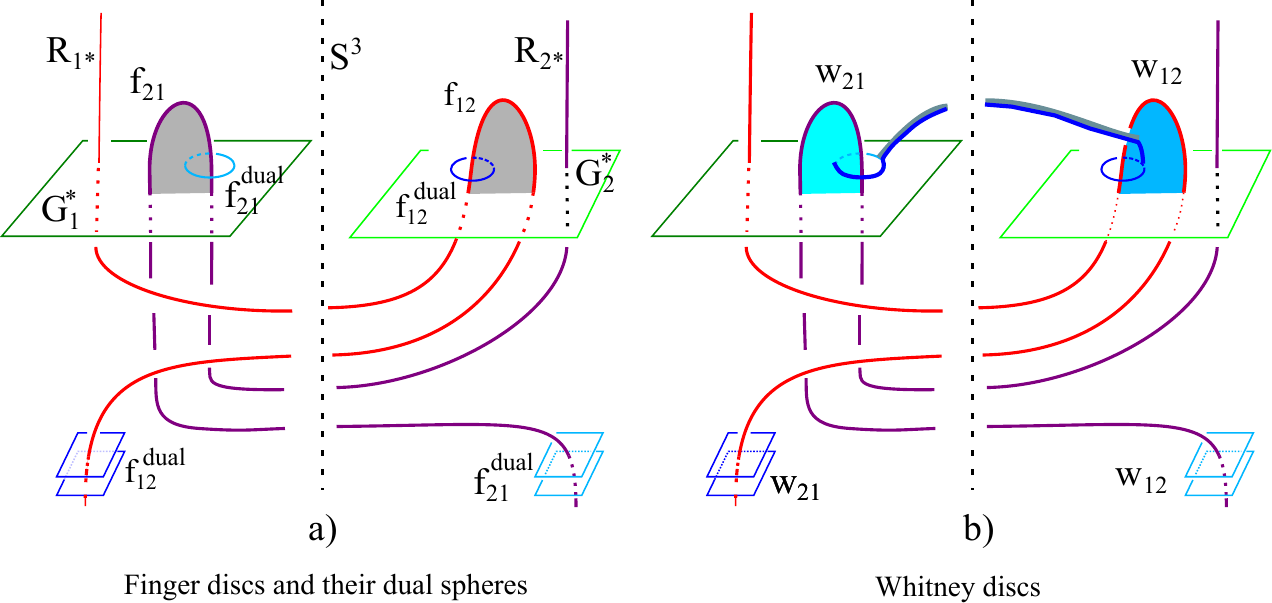}
    \caption{The double cross family viewed from the middle middle level}
    \label{fig:doublecross1}
\end{figure}

\vskip 8pt
In this section we define the \emph{double cross family}, a 2-eyed Hatcher-Wagoner family with graphic as in Figure \ref{fig:graphic} and finger-Whitney system $(\mF_D,\mW_D)$  shown in Figures \ref{fig:doublecross1} a), b),    discovered jointly by the authors and David Gay.    We will first explain how to decipher Figure \ref{fig:doublecross1}, then explicitly describe a 2-parameter deformation which transforms this family to one with finger-Whitney system $\mfw$ that of the barbell family.   This shows that the double cross family induces the barbell diffeomorphism of the 4-sphere.

\vskip 8pt
In \cite{GGH}, two functions $\I$ and $\C$ are defined on finger-Whitney systems for loops of spheres in $\#^k S^2\times S^2$ and each of these functions count certain intersections between finger disks and Whitney disks $\mod 2$. With respect to the functions $C$ and $\I$ on finger-Whitney systems, $C\mfw=0$ and $\I\mfw=1$, see \S0 \cite{GGH1}.   By Proposition \ref{d double cross},  $C(\mF_D,\mW_D)=1$ and $\I(\mF_D, \mW_D)=0$.  Thus, the quantities $\I$ and $C$ may change under 2-parameter deformation.  However, \cite{GGH1} shows that $\mD=C+\I$ is invariant under deformations through Hatcher-Wagoner families and the sequel \cite{GGH2} shows that $\mD$ is an invariant of Hatcher-Wagoner families connected by general 2-parameter families.

\vskip 8pt
We now explain how to read off the finger-Whitney system depicted in Figures \ref{fig:doublecross1} a) and \ref{fig:doublecross1} b), the finger disks and their dual spheres being shown in Figure \ref{fig:doublecross1} a) and the Whitney disks in Figure \ref{fig:doublecross1} b).  These figures show a 3D slice of the middle middle level $V_x=q^{-1}_{\frac{1}{2}}(\frac{1}{2}) =\stwostwo\#\stwostwo$, where $x=(1/2,1/2)\in$ the graphic and $q_{\frac{1}{2}}:S^4\times I\to I$ is the corresponding Morse function.  Here $G_1, G_2$ are the ascending spheres of the critical points of index $2$ denoted $G_1, G_2$ in the graphic.  Similarly the spheres $R_{1}, R_{2}$ are the descending spheres of the index $3$ critical points.  We use the notation $f_{ij} $ to indicate that it is a finger disk between $R_{i} $ and $G_j$. 
    
We  describe the dual spheres schematically shown in Figure \ref{fig:doublecross1} a).  Here $f_{21}^{\textrm{dual}}$ is built using Quinn's construction \cite{Qu}.  It satisfies $|f_{21}^{\textrm{dual}}\cap f_{21}|=1,  f_{21}^{\textrm{dual}}\cap(R_{1}\cup R_{2}\cup G_1\cup G_2\cup f_{12}\cup f_{12}^{\textrm{dual}})=\emptyset$ and    $f_{12}^{\textrm{dual}}$ is constructed similarly with corresponding properties.  In more detail to construct $f_{21}^{\textrm{dual}}$, start with a round 2-sphere which intersects the 3D slice in the circle intersecting $f_{21}$ as shown in Figure \ref{fig:doublecross1} a).  That 2-sphere intersects $R_{2}$ once in the future and once in the past.  Each intersection is tubed off with a dual sphere to $R_{2}$, both of which are schematically shown in Figure \ref{fig:doublecross1} a).  

Next we describe the Whitney disks $w_{12}$ and $w_{21}$, where $w_{ij}$ denotes the Whitney disk between $R_{i}$ and $G_j$.  Construct $w_{12}$ by tubing a copy of $f_{12}$ to $f_{21}^{\textrm{dual}}$.  Note that $|w_{12}\cap f_{21}|=1$ and $|w_{12}\cap f_{12}^{\textrm{dual}}|=1$.  Next construct $w_{21}$ by tubing a copy of $f_{21}$ to $f_{12}^{\textrm{dual}}$ by a tube that follows an arc in $w_{12}$ from $w_{12}\cap f_{21}$ to $w_{12}\cap f_{12}^{\textrm{dual}}$.

\begin{proposition} \label{d double cross}  $\I(\mF_D,\mW_D)=0$ and $C(\mF_D,\mW_D)=1$.  \end{proposition}
\begin{proof}  After a suitable preliminary isotopy of $\mW_D$ to be Whitney germed and $\mF_D$ to be finger germed, $(\mF_D,\mW_D)$ is in full embedded arc position.  Since $(\mF_D,\mW_D)$ has no uncross disks and for all $i\neq j, |f_{ij}\cap w_{ji}|=1$ the conclusion follows from \cite{GGH1}.  \end{proof}

\begin{theorem}  \label{double cross barbell}  The double cross family induces the barbell diffeomorphism, indeed there is a 2-parameter family of Morse function and gradient like vector fields interpolating between the double cross family of spheres and the barbell family of spheres.\end{theorem}

\begin{proof}  View the double cross family from the $y$-middle level $V_y=\stwostwo$.  See Figures \ref{fig:doublecross2} a) and \ref{fig:doublecross2} b).  Here $V_y$ is the level surface of $\pi_{\frac{1}{2}}$ above the first index $3$ critical point, where $y$ is seen in the graphic.  Also $V_y$ and the data of   Figures \ref{fig:doublecross2} a) and \ref{fig:doublecross2} b) are respectively obtained from $V_x$, Figures \ref{fig:doublecross1} a) and \ref{fig:doublecross1} b) by surgering $R_{2*}$.  

\begin{figure}
    \centering
    \includegraphics[width=1.\linewidth]{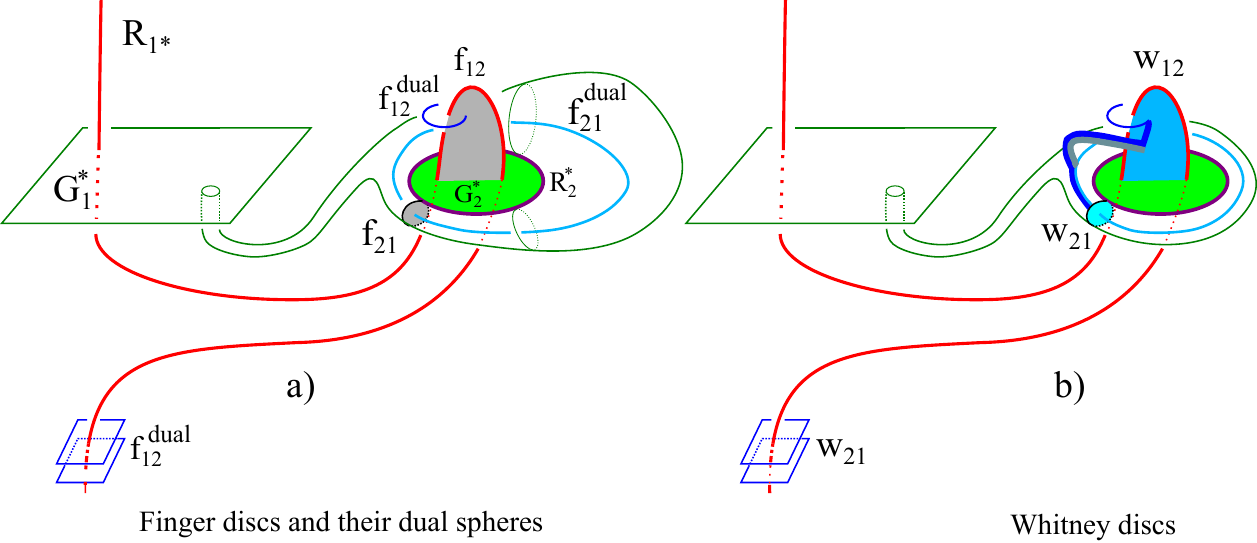}
    \caption{The double cross family viewed from the $y$-middle level}
    \label{fig:doublecross2}
\end{figure}

\begin{figure}
    \centering
    \includegraphics[width=.9\linewidth]{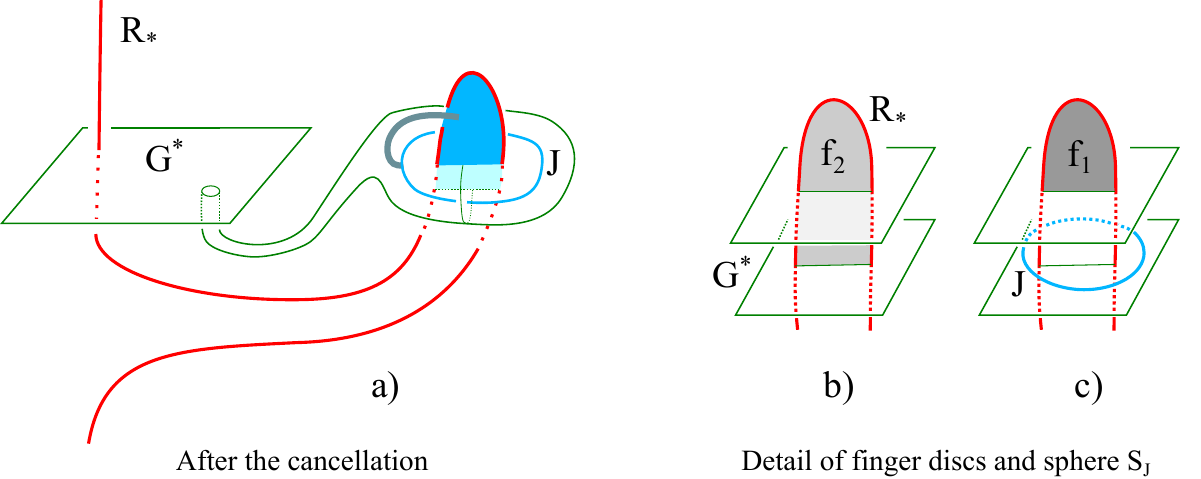}
    \caption{Reduction to one eye}
    \label{fig:one eye reduction}
\end{figure}
\vskip 8pt
We next cancel the inner eye.  This can be done since  there is a unique flow line from $R_2$ to $G_2$ throughout the eye.  The result of the cancellation is a new family with a single eye.  Figure \ref{fig:one eye reduction} depicts its middle middle level with $G$ (resp. $R$) the ascending (resp. descending) sphere of the critical point of index $2$ (resp. index $3$).  This family has two finger disks denoted by $f_1$ and $f_2$ respectively shown in Figures \ref{fig:one eye reduction} c) and \ref{fig:one eye reduction} b).  We can assume that $f_2$ is  slightly into the future to be disjoint from $f_1$.
\begin{figure}
    \centering
    \includegraphics[width=.7 \linewidth]{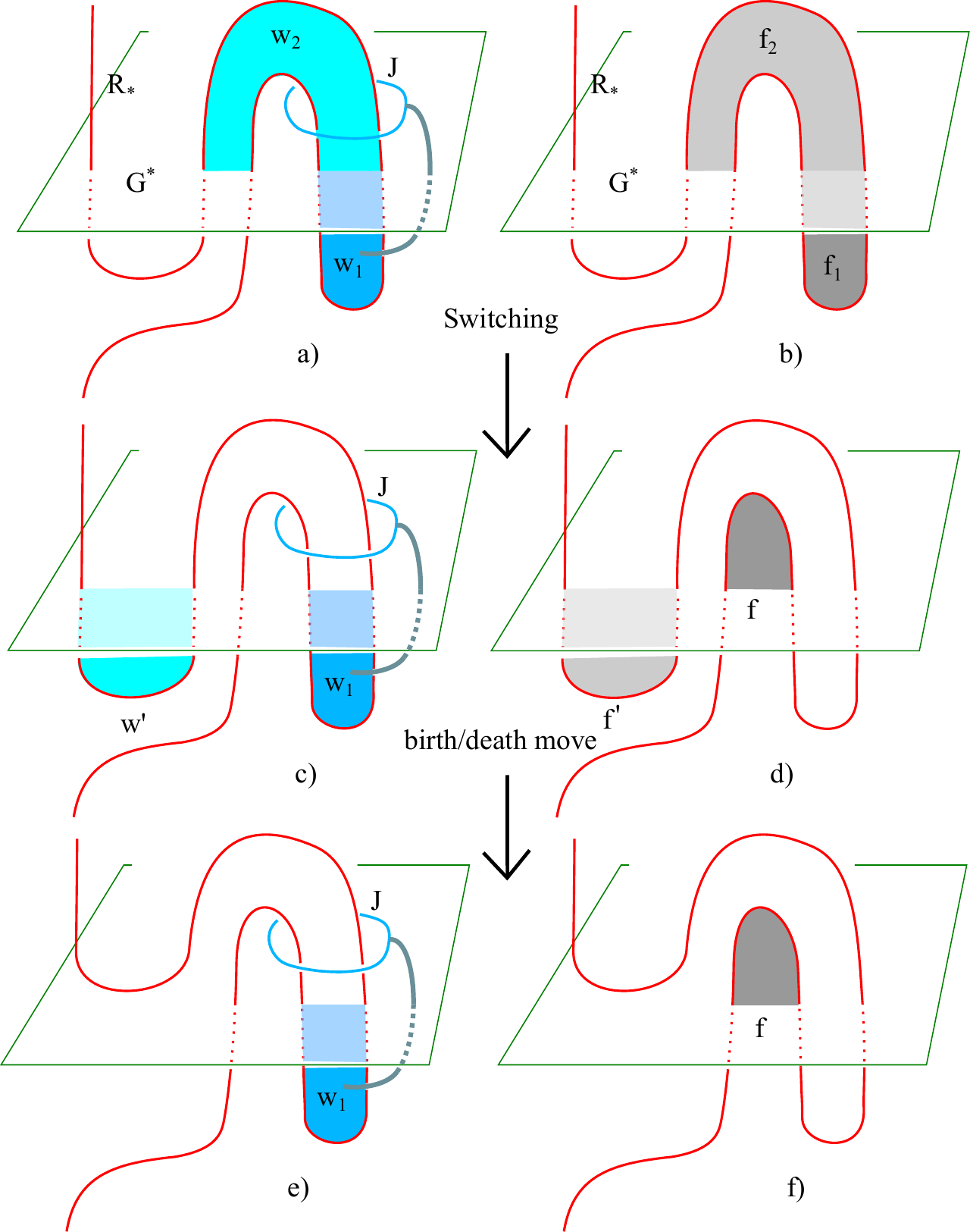}
    \caption{Reduction to the barbell family}
    \label{fig:barbell reduction}
\end{figure}
Here $J$ is the intersection with the 3D slice of a 2-sphere $S_J$ disjoint from $R*\cup G^*\cup f_1\cup f_2$.  It bounds disks in the future and past.  The Whitney disks, denoted $w_1$ and $w_2$ are obtained respectively from $f_1$ and $f_2$ by tubing to parallel copies of $S_J$ as in Figure \ref{fig:one eye reduction} a) using parallel tubes.  

Figure \ref{fig:barbell reduction} is the result of isotoping Figure \ref{fig:one eye reduction}, where Figure \ref{fig:barbell reduction} a) (resp. \ref{fig:barbell reduction} b)) shows the Whitney (resp. finger) disks.  Figure \ref{fig:barbell reduction} c) is obtained from Figure \ref{fig:barbell reduction} a) by switching $w_2$ to $w'$ and Figure \ref{fig:barbell reduction} d) is obtained from \ref {fig:barbell reduction} b)  by switching $f_2$ to $f'$ and then $f_1$ to $f$.  These switchings are 1-switches in the terminology of \S10 \cite{GGH1}.  The are realizable by a 2-parameter family by Lemma 0.55 \cite{GGH1}.  Since $w'=f'$ we can apply the cancellation move of \cite{GGH1} to eliminate them thereby obtaining a finger-Whitney system whose Whitney and finger disks are respectively shown  in Figures \ref {fig:barbell reduction} e) and \ref{fig:barbell reduction} f).  Since this is a finger-Whitney system for the barbell family the result follows.\end{proof}

\newpage

\printbibliography[heading=bibintoc]

\end{document}